\documentclass[11pt]{amsart}
\usepackage{latexsym}
\usepackage{graphicx}
\begin{document}

\newtheorem{thm}{Theorem}[section] 
\newtheorem{cor}[thm]{Corollary}
\newtheorem{lem}[thm]{Lemma}
\newtheorem{prop}[thm]{Proposition}
\newtheorem{rem}[thm]{Remark}
\newtheorem{defn}[thm]{Definition}
\newtheorem{defs}[thm]{Definitions}
\newtheorem{con}{Conjecture}
\newcommand{\aaa}{\mbox{$\alpha$}}
\newcommand{\map}{\mbox{$\rightarrow$}}
\newcommand{\up}{\mbox{$\uparrow$}}
\newcommand{\down}{\mbox{$\downarrow$}}
\newcommand{\ccc}{\mbox{$\mathcal C$}}
\newcommand{\kkk}{\mbox{$\kappa$}}
\newcommand{\Aaa}{\mbox{$\mathcal A$}}
\newcommand{\Baa}{\mbox{$\mathcal B$}}
\newcommand{\Daa}{\mbox{$\mathcal D$}}
\newcommand{\Fff}{\mbox{$\mathcal F$}}  
\newcommand{\bbb}{\mbox{$\beta$}}
\newcommand{\sss}{\mbox{$\sigma$}}  
\newcommand{\ddd}{\mbox{$\delta$}} 
\newcommand{\rrr}{\mbox{$\rho$}} 
\newcommand{\ooo}{\mbox{$\omega$}} 
\newcommand{\Ggg}{\mbox{$\Gamma$}}
\newcommand{\gggg}{\mbox{$\gamma$}}
\newcommand{\ttt}{\mbox{$\tau$}} 
\newcommand{\bdd}{\mbox{$\partial$}}
\newcommand{\zzz}{\mbox{$\zeta$}}
\newcommand{\Ss}{\mbox{$\Sigma$}}
\newcommand{\Sss}{\mbox{$\mathcal S$}}
\newcommand{\Ddd}{\mbox{$\Delta$}}
\newcommand{\Fa}{\mbox{$\Fff_A$}}
\newcommand{\Fb}{\mbox{$\Fff_B$}}
\newcommand{\llll}{\mbox{$\lambda$}}
\newcommand{\pf}{\noindent{\bf Proof: }}
\renewcommand{\baselinestretch}{1.5}


\title[Annuli and Heegaard splittings] {Annuli in generalized Heegaard
splittings and degeneration of tunnel number}

\author{Martin Scharlemann}
\address{\hskip-\parindent
        Mathematics Department\\
        University of California\\
        Santa Barbara, CA 93106}
\email{mgscharl@math.ucsb.edu}

\author{Jennifer Schultens}
\address{\hskip-\parindent
        Department of Mathematics \& CS\\
        Emory University \\
        Atlanta, GA 30322}
\email{jcs@@mathcs.emory.edu}

\thanks{The authors are supported in part by National Science Foundation
grants.}

\date{\today}

\begin{abstract} We analyze how a family of essential
annuli in a compact $3$-manifold will induce, from a strongly irreducible
generalized Heegaard splitting of the ambient manifold, generalized
Heegaard splittings of the complementary components.  There are specific
applications to the subadditivity of tunnel number of knots, improving
somewhat bounds of Kowng \cite{Kw}.  For example, in the absence of
$2$-bridge summands, the tunnel number of the sum of $n$ knots is no less
than $\frac{2}{5}$ the sum of the tunnel numbers.
\end{abstract}

\maketitle

\section{Introduction}

\vspace{5 mm}

The {\em tunnel number} $t(K)$ of a knot (or link) $K$ in $S^3$ is the
minimal number of arcs that, when attached to the knot, gives a graph
whose complement is the interior of a handlebody.  It was once naively
hoped that this knot invariant might be additive under connected sum of
knots, but Morimoto (\cite{M2}) has found counterexamples in which tunnel
number  {\em degenerates} under
connected sum.  That is, there are knots $K^1$ and $K^2$
for which $t(K^1 \# K^2) < t(K^1) + t(K^2)$.  The degree of degeneration
possible was soon shown by Kobayashi (\cite{Ko}) to be arbitrarily high. 
That is, given $d \geq 0$ there exist knots $K^1$ and $K^2$ for which
$t(K^1) + t(K^2) - t(K^1 \# K^2) > d$. 

Although Kobayashi's examples show that there is no uniform
bound on the degeneration number, it is still natural to ask if
there is a uniform bound less than $1$ on the {\em degeneration ratio}
$$d(K^1, K^2) = \frac{t(K^1) + t(K^2) - t(K^1 \# K^2)}{t(K^1) +
t(K^2)}.$$  More generally, for knot summands $K^1, \ldots , K^n$ we can
define the degeneration ratio $$d(K^1, \ldots , K^n) = \frac{t(K^1) +
\ldots + t(K^n) - t(K^1 \#  \ldots \# K^n)}{t(K^1) + \ldots +
t(K^n)}.$$  Kobayashi's methods show that the degeneration ratio for each
of his examples is at least
$\frac{1}{9}$ and Morimoto's specific example has degeneration ratio
$\frac{1}{3}$.  If we also consider links, Morimoto has shown
\cite{M1} that  the connected sum of a knot $K$ and a $2$-component link
$L$ has tunnel number one exactly if $K$ is a $2$-bridge knot and $L$ is
the Hopf link.  Such examples have degeneration ratio $\frac{1}{2}$.

In the other direction, Kowng showed in \cite{Kw} the more general claim
that if an irreducible orientable compact $3$-manifold $M$ contains a
collection of tori
$\mathcal T$, then $M$ has a generalized Heegaard splitting of genus at most
$3(genus(M) + | {\mathcal T} |) -2$ for which each component of
$\mathcal T$ is contained in a thin level (see Definition \ref{general}). 
When applied to $M = S^3 - \eta(K^1 \#  \ldots \# K^n)$, this
has the corollary that $$\frac{1}{3} (t(K^1) + \ldots +t(K^n)) \leq t(K^1
\#  \ldots \# K^n) + (n-1).$$  So Kowng showed that, in general, 
$$d(K^1, \ldots , K^n) \leq \frac{2}{3} + \frac{n-1}{t(K^1) + \ldots
+t(K^n)}.$$

Here we improve this aspect of Kowng's results, using a somewhat different
approach growing out of the investigations in \cite{SS1}, \cite{SS2}. 
The Main Theorem shows how to cut a strongly irreducible generalized
Heegaard splitting of $M$ along a family of essential annuli and create
generalized Heegaard splittings for the resulting manifolds. Analysis of
the result shows a bit more than the following: for $K^j$ prime, $d(K^1,
\ldots , K^n) \leq \frac{2}{3}$; if none of the $K^j$ is $2$-bridge, then
$d(K^1, \ldots , K^n) \leq \frac{3}{5}$.

$2$-bridge knots seem to play a special role in the theory of tunnel
numbers.  Morimoto's original example of degeneration used a $2$-bridge
summand, and Kobayashi's (non-prime) examples each have for one of their
summands a sum of $2$-bridge knots.  So it is intriguing that here too
the bounds on degeneration are better if there are no $2$-bridge knots
among the summands. 

Although we focus on cutting a manifold apart by annuli, the same sort of
analysis extends to families of tori, in a way we only briefly
describe.

\bigskip

After some preliminaries and a simplified but indicative special case, the
outline of the proof is as follows:

{\bf Step 1 (Section \ref{sec:annuli}):}  We examine annuli in, and on the
boundary of, compression bodies.  The goal is to develop criteria that
allow us to cut open compression bodies along annuli and glue compression
bodies together along annuli so that at the end, the result is still a
collection of compression bodies.

{\bf Step 2 (Section \ref{sec:break}):} We prove that a
strongly irreducible generalized Heegaard splitting can be cut along
annuli and augmented in such a way as to yield generalized Heegaard
splittings for the complementary components.  In this construction a
bound on the degeneration of the ``index'' is given by the number of
``dipping annuli''.

{\bf Step 3 (Section \ref{sec:des}):} We show that the generalized
Heegaard splittings of the complementary components obtained in this
way can be destabilized under certain conditions.

{\bf Step 4 (Sections \ref{sec:generalcount} and \ref{sec:specialcount}):}
We find bounds on the number of dipping annuli.  Section
\ref{sec:generalcount} is quite general, whereas Section
\ref{sec:specialcount} achieves better bounds in the context of a knot
complement, when we can assume any surface with meridinal boundary has an
even number of boundary components.

{\bf Step 5 (Section \ref{sec:tun}):} We apply the derived inequalities
to the study of tunnel numbers.

\bigskip

The paper concludes with an appendix by Andrew Casson.  He constructs an
example which demonstrates that the bound on the number of dipping annuli
given in Section 7 is, in some sense, best possible.

\section{Preliminaries}

For standard definitions concerning knots, see \cite{BZ} or \cite{R} and
for those concerning $3$-manifolds, see \cite{H} or \cite{J}.  All
manifolds will be orientable.

\begin{defn} For $N$ a properly embedded submanifold of $M$, we denote an
open regular neighborhood of $N$ in $M$ by $\eta(N)$.
\end{defn}

\begin{defn} Let $K$ be a knot in $S^3$.  Denote the
\underline{complement of K},
$S^3 - \eta(K)$, by $C(K)$.  
\end{defn}

\begin{rem}  \label{rem:dec} Let $K = K^1 \# K^2$ be the sum of two
knots.  Then the decomposing sphere gives rise to a decomposing annulus
$A$ properly embedded in $C(K)$ such that $C(K) = C(K^1) \cup_A C(K^2)$. 
If $K = K^1 \# \dots \# K^n$, then we may assume that the decomposing
spheres are nested, so that $C(K) = C(K^1) \cup_{A_1} \dots
\cup_{A_{n-1}} C(K^n)$.  
\end{rem}

\begin{defn} \label{defn:tunnel system} A \underline{tunnel system} for a
knot $K$ is a collection of disjoint arcs ${\mathcal T} = t_1 \cup \dots \cup
t_n$ properly embedded in $C(K)$ such that $C(K) - \eta({\mathcal T})$ is a
handlebody.  The
\underline{tunnel number of $K$}, denoted by $t(K)$, is the least number
of arcs required in a tunnel system for $K$.   
\end{defn}
 
\begin{defn} \label{defn:cb} A \underline{compression body} is a
$3$-manifold $W$ obtained from a connected closed orientable surface $S$
by attaching $2$-handles to
$S \times \{0\} \subset S \times I$ and capping off any resulting
$2$-sphere boundary components.  We denote $S \times \{1\}$ by
$\partial_+W$ and $\partial W - \partial_+W$ by
$\partial_-W$.  Dually, a compression body is a connected
orientable $3$-manifold obtained from a (not necessarily connected) closed
orientable surface
$\partial_-W \times I$ by attaching $1$-handles.  Define the
\underline{index} of $W$ by $J(W) = \chi(\partial_-W) -
\chi(\partial_+W) \geq 0$.
\end{defn}

\begin{defn} A \underline{set of defining disks} for a compression body
$W$ is a set of disks $\{D_1, \dots, D_n\}$ properly embedded in $W$ with
$\partial D_i \subset \partial_+W$ for $i = 1$, $\dots, n$ such that the
result of cutting $W$ along $D_1 \cup \dots \cup D_n$ is homeomorphic to
$\partial_-W \times I$.
\end{defn}

\begin{defn} \label{defn:Heegaard splitting} A \underline{Heegaard
splitting} of a $3$-manifold $M$ is a decomposition $M = V \cup_S W$ in
which $V$, $W$ are compression bodies such that $V \cap W = \partial_+V =
\partial_+W = S$ and $M = V \cup W$.  We call $S$ the
\underline{splitting surface} or \underline{Heegaard surface}.
\end{defn}

\begin{defn} \label{defn:wr} A Heegaard splitting is
\underline{reducible} (resp. \underline{weakly} \underline{reducible}) if
there are essential disks $D_1$ and $D_2$,  such that
$\partial D_1 = \partial D_2$ (resp. $\partial D_1 \cap \partial D_2 =
\emptyset$).  A Heegaard splitting which is not (weakly) reducible
is \underline{(strongly)} \underline{irreducible}. \end{defn}

\begin{defn}
\label{general} A \underline{generalized} Heegaard splitting of a compact
orientable $3$-manifold $M$ is a structure $M = (V_1 \cup_{S_1} W_1)
\cup_{F_1} (V_2 \cup_{S_2} W_2) \cup_{F_2} \dots \cup_{F_{m-1}} (V_m
\cup_{S_m} W_m)$.  Each of the $V_i$ and
$W_i$ is a union of compression bodies, $\partial_+V_i = S_i =
\partial_+W_i$, (i.e., $V_i \cup_{S_i} W_i$ is a union of Heegaard
splittings of a submanifold of
$M$) and $\partial_-W_i = F_i = \partial_-V_{i+1}$.   We say that a
generalized Heegaard splitting is
\underline{strongly} \underline{irreducible} if each Heegaard splitting
$V_i
\cup_{S_i} W_i$ is strongly irreducible and each $F_i$ is
incompressible in $M$.  We will denote $\cup_i F_i$ by
$\mathcal F$ and
$\cup_i S_i$ by $\mathcal S$.  The surfaces in $\mathcal F$ are called the
\underline{thin levels} and the surfaces in $\mathcal S$ the
\underline{thick levels}.

Let $M = V \cup_S W$ be an irreducible Heegaard splitting.  We may think
of $M$ as being obtained from $\partial_-V \times I$ by attaching all
$1$-handles in $V$ followed by all $2$-handles in $W$, followed, perhaps,
by $3$-handles.  An \underline{untelescoping} of
$M = V \cup_S W$ is a rearrangement of the order in which the
$1$-handles of $V$ and the $2$-handles of $W$ are attached.  This
rearrangement yields a generalized Heegaard splitting.  For convenience,
we will occasionally denote $\partial_-V = \partial_-V_1$ by $F_0$.
\end{defn}

The Main Theorem in \cite{ST1} together with the calculation
\cite[Lemma 2]{SS2} implies the following:

\begin{thm} \label{lem:ind} Suppose $M$ is an irreducible compact
$3$-manifold.  Then $M$ possesses a genus $g$ Heegaard
splitting if and only if $M$ has a strongly irreducible generalized
Heegaard splitting
$(V_1 \cup_{S_1} W_1) \cup_{F_1} (V_2 \cup_{S_2} W_2) \cup_{F_2} \dots
\cup_{F_{m-1}} (V_m \cup_{S_m} W_m)$ such that $\sum_{i=1}^m J(V_i) =
\sum_{i=1}^m J(W_i) = 2g - 2$.
\end{thm}

One implication comes from untelescoping a Heegaard splitting of genus
$g$, the other from thinking of a generalized Heegaard splitting as an
untelescoping of some Heegaard splitting.  The latter process is
sometimes called the {\em amalgamation} of the generalized splitting into
a standard splitting.

A strongly irreducible Heegaard splitting can be isotoped so that its
splitting surface, $S$ intersects an incompressible surface, $P$, only in
curves essential in both $S$ and $P$.  This is a deep fact and is proven,
for instance, in  \cite[Lemma 6]{Sc}.  This fact, together with the fact
that incompressible surfaces can be isotoped to meet only in essential
curves, establishes the following:

\begin{lem}  \label{lem:essential} Let $P$ be a properly embedded
incompressible surface in an irreducible $3$-manifold $M$ and let $M =
(V_1 \cup_{S_1} W_1) \cup_{F_1} \dots \cup_{F_{m-1}} (V_m \cup_{S_m}
W_m)$ be a strongly irreducible generalized Heegaard splitting of
$M$.  Then ${\mathcal F} \cup {\mathcal S}$ can be isotoped to intersect
$P$ only in curves that are essential in both $P$ and $\Fff \cup \Sss$.
\end{lem}

\section{An indicative first example}

We will first point out a special simplified case in which an upper bound
to the degeneration ratio is fairly easily found, and the ideas that are
used are indicative of ideas that will be important in the general case.

The special case is this:

\begin{thm} \label{thm:indic}Let $K^1$ and $K^2$ be prime knots and assume
that
$M = C(K^1
\# K^2)$ possesses a minimal genus Heegaard splitting
that is strongly irreducible.  Then $t(K^1 \# K^2) \geq \frac{1}{2}(t(K^1)
+ t(K^2))$.  Hence $d(K^1, K^2) \leq \frac{1}{2}$.
\end{thm}

\pf
Let $M = V \cup_S W$ be the strongly irreducible Heegaard splitting,
with $\bdd M = \bdd_- V$ and $W$ a handlebody.  The hypothesis implies 
that
$2t(K^1
\# K^2) = -
\chi(S)$.  Let
$(A,
\bdd A) \subset (M, \bdd M)$ denote the annulus in $M = C(K^1 \# K^2)$
that, when completed by a pair of meridian disks, constitutes a
decomposing sphere for $K^1 \# K^2
\subset S^3$.  Following Lemma \ref{lem:essential} we can isotope $S$ and
$A$ until their intersection consists of circles essential in $A$ and in
$S$.  Isotope them further to reduce, as much as possible, the number of
such components of $A \cap S$. 

Let $A_V = A \cap V$ and $A_W = A \cap W$ be the essential annuli $A - S$
in $V$ and $W$ respectively.  A sequence of $\bdd$-compressions of $A_W
\subset W$ could turn each annulus into a disk; similarly a sequence of
$\bdd$-compressions of all but the spanning annulus of $A_V \subset V$
could turn each annulus into a disk.  If these sets of $\bdd$-compressions
could be done simultaneously, the result would be that $A \cap S$ would
become a single circle of intersection, inessential in $A$ and bounding
an essential disk in $W$.  The disk would divide $W$ into two
handlebodies, and it would follow that $t(K_1) + t(K_2) \leq t(K^1 \#
K^2) - 1$ so there would be no degeneration (in fact, a surplus!)

In general, the boundary compressions of $A_V$ and $A_W$ cannot be done
simultaneously, since the arcs to which they $\bdd$-compress in $S$ may
intersect.  This problem can be avoided if
$S$ were first stabilized by attaching a tube parallel to a spanning arc
of each annulus in $A_V$.  This would allow
$\bdd$-compressions of $A_V$ onto arcs that run along the tubes, without
affecting the $\bdd$-compressions of $A_W$ to the original $S$.  So we
see immediately that 
$|A_V|$, the number of tubes used, is an upper bound on degeneration. The
first goal is then to get a bound on $|A_V|$.  

Since any surface in $S^3$ is separating, any component of $S -
A$ has an even number of boundary components, since it can be
completed to become a closed surface in $S^3$ by attaching meridian disks
of either $K^1$ or $K^2$.   Assume, for initial simplification, that no
component of $S - A$ is an annulus.  Let $S' = M - \eta(A)$ and notice
then that each component $S_0$ of $S'$ has at most $2 - \chi(S_0) \leq -2
\chi(S_0)$ boundary components, since each component has non-trivial
even Euler characteristic.  It follows that $|\bdd S'| \leq -2 \chi(S') =
-2
\chi(S)$, so the number of essential curves $|A \cap S| \leq
-\chi(S)$.  This implies that $|A_V| \leq -\frac{\chi(S)}{2} = t(K^1
\# K^2)$.  Combining this argument with the previous one, we get $$t(K_1)
+ t(K_2) \leq  t(K^1 \# K^2) + |A_V| - 1 \leq 2 t(K^1 \# K^2) - 1$$ or
$$t(K^1 \# K^2) > \frac{1}{2}(t(K_1) + t(K_2)),$$ an inequality better
than required.

Now examine the annulus components in $S - A$.  Note first of all that
any such annulus component $B$ is necessarily $\bdd$-parallel in the
component $C(K^1)$ (say) of $M - A$ in which it lies, since 
$K^1$ is prime.  Since the number of (essential) circles $|S \cap A|$ has
been minimized, $B$ is not parallel to a subannulus of $A$.  So $B$ must
be parallel to the annulus component of $\bdd M - \bdd A$ that lies
in $C(K^1)$, an annulus on $\bdd M$ which we denote by $B_{\bdd}$.  

It
follows that there cannot be two annulus components of $S - A$ that are
adjacent in $S$, since one would be parallel to each of the two annulus
components of $\bdd M - \bdd A$ and so the union of the two annuli in
$S$ would provide a way to ''spin'' a collar of $\bdd A$ in a way that
would reduce the number of components of $S \cap A$.  Since no two annuli
in $S - A$ are adjacent in $S$, the total number of these annuli can be no
larger than the number of circles $|A \cap S|$ that we calculated above
in the absence of annuli.  That is, there are no more than
$-\chi(S)$ annuli among the components of $S - A$.  Hence, even allowing
annuli components in $S - A$,
$|A
\cap S| \leq -2\chi(S)$ so $|A_V| \leq \chi(S) = 2t(K^1 \# K^2)$ and we
get $$t(K_1) + t(K_2) \leq  t(K^1 \# K^2) + |A_V| - 1 \leq 3 t(K^1 \# K^2)
- 1$$ or $$t(K^1 \# K^2) > \frac{1}{3}(t(K_1) + t(K_2)).$$ 

But we can do better.  The region lying between the annuli $B$
and $B_{\bdd}$ in $C(K^1)$ is homeomorphic to $annulus \times I$, i. e. a
solid torus, and it is known how a strongly irreducible splitting surface
like
$S$ can intersect a solid torus (cf \cite{Sh}).  The upshot is that all
components of $S - A$ lying between $B$ and $B_{\bdd}$ in $C(K^1)$ are
annuli parallel to $B$, except possibly one component (an {\em
exceptional} component) which consists of a pair of annuli parallel to $B$
but then tubed together by a vertical tube.  (See Figure 1.)  An
exceptional component can be ignored since, for example, the vertical
tube could be slid across $A$ into
$C(K^2)$.  Notice that, by choosing $B$ outermost (i. e. furthest from
$B_{\bdd} \subset \bdd M$), all annuli components of $S \cap C(K^1)$ lie
between
$B$ and
$B_{\bdd}$.  By the analogous argument in
$C(K^2)$ no component of $S \cap C(K^2)$ can be an annulus, since an
innermost one would be adjacent in $S$ to an innermost one of $S \cap
C(K^1)$, forming a $\bdd$-parallel torus component of $S$, a
contradiction.  (This argument is complicated only a little by the
possible presence of an exceptional component.)

\begin{figure}
\centering
\includegraphics[width=.8\textwidth]{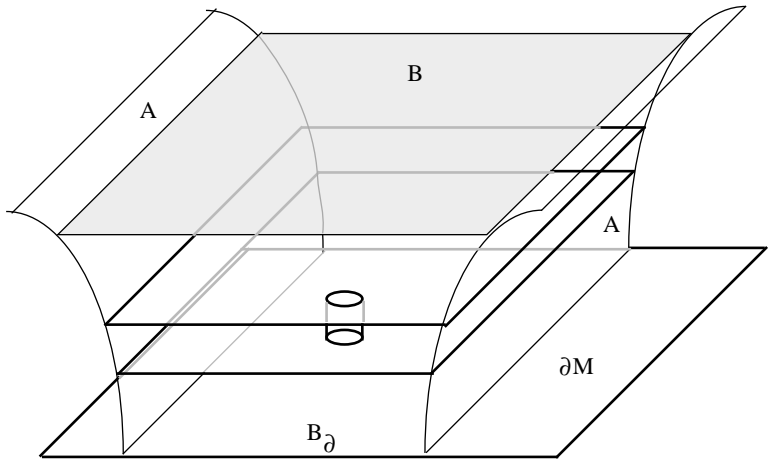}
\caption{An exceptional component between $B$ and $B_{\bdd}$}
\end{figure} 

So we may suppose that all $w$ annulus components of $S - A$ lie in
$C(K^1)$, lie between and are parallel to the annuli $B$ and $B_{\bdd}$.
Then at least
$\frac{w-2}{2}$ pairs of these annuli each cut off a parallel pair of
annuli in $A_V$.  For each of these pairs, only one stabilization is
required in the process described above, since the annuli are parallel.
(Imagine placing the added tube so that it lies in the collar between the
annuli in $V$). So if $S - A$ contains $w$ annuli there are two new
effects:  First, the count of
$|\bdd S'|$ is changed to 
$|\bdd S'| \leq 2w -2\chi(S') = 2w -2 \chi(S)$, so the number of essential
curves $|A
\cap S| \leq w-\chi(S)$.  This implies that $$|A_V| \leq \frac{w -
\chi(S)}{2}.$$  On the other hand, the second effect is that the number of
stabilizations required is at most $$|A_V| - \frac{w-2}{2} \leq
\frac{-\chi(S)}{2} + 1 = t(K^1 \# K^2) + 1.$$  This implies that $$t(K_1)
+ t(K_2) \leq  t(K^1
\# K^2) + (t(K^1 \# K^2) + 1) - 1 \leq 2 t(K^1 \# K^2)$$ or
$$t(K^1 \# K^2) \geq \frac{1}{2}(t(K_1) + t(K_2)),$$ as required.
\qed

\section{Annuli in and on compression bodies} \label{sec:annuli}

A properly imbedded annulus $A \subset M$ is {\em essential} if it is
incompressible and not $\bdd$-parallel.  In this section we study
finite sets of disjoint essential annuli in a compression body $W$.

\begin{lem}  \label{lem:br} If $A$ is an essential annulus in a
compression body $W$, then either one component of $\bdd A$ lies on
each of $\bdd_+ W$ and $\bdd_- W$ (i. e. $A$ is {\em spanning}) or 
$\bdd A \subset \bdd_+ W$.  In the latter case, $A$ is boundary
compressible.  
\end{lem}

\pf See \cite[Lemma 9]{BnO}. \qed
\bigskip

Suppose \Aaa\ is a properly imbedded collection of essential annuli in a
compression body $W$ so, in particular, the boundary of any non-spanning
annulus $A$ lies in $\bdd_+ W$.  Suppose $D$ is a
$\bdd$-compressing disk for the non-spanning annulus $A$, i. e. the
interior of $D$ is disjoint from $\Aaa$, and $\bdd D$ is the union of a
spanning arc \aaa\ of $A$ and an arc in $\bdd_+ W$. The arc components of
$D \cap \Aaa$ divide $D$ into subdisks.  This naturally gives rise to a
tree in $D$, in which a vertex is chosen inside each subdisk and two
vertices are connected if they abut the same annulus.  (We can ignore
closed components of intersection.  Indeed, since $\Aaa$ is
incompressible,  closed components of  $D \cap \Aaa$ can be removed by an
isotopy with support disjoint from
the arcs of intersection.) It will be useful to extend this tree by
attaching, to the vertex corresponding to the subdisk that abuts \aaa, an
edge \rrr\ that crosses \aaa.  The other end of \rrr\ is called the root
of the resulting tree \ttt; the other valence one edges are called the
leaves of \ttt.  Each leaf corresponds to a disk cut off by an outermost
arc of
$\Aaa \cap D$ in $D$. (See Figure 2.)

\begin{figure}
\centering
\includegraphics[width=.8\textwidth]{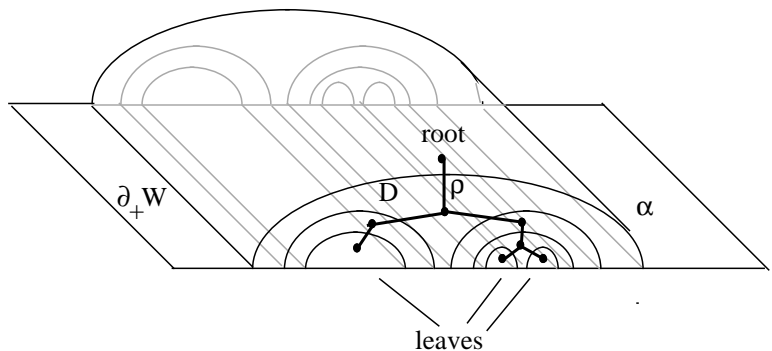}
\caption{A tree of complexity $\{ 4,4,3 \}$.}
\end{figure} 

Recall that there is a natural order on the set of finite sets of
integers (see \cite[Definition 4.3]{Ga}).  One set of
integers is compared to another by arranging each in descending (or at
least never-ascending) order and then comparing them lexicographically. 
This ordering has the property that if a subset of a set of integers is
replaced by a subset of lower order then the resulting set is of lower
order.  

\begin{defn}   The \underline{complexity} of the tree \ttt\ is the set of 
distances (measured in edges traversed) from the root of \ttt\ to the
leaves in \ttt.  The complexity of a $\bdd$-compressing disk for $A \in
\Aaa$ is the complexity of the associated tree.  A \underline{minimal}
$\bdd$-compressing disk for $A$ is a $\bdd$-compressing disk of minimal
complexity (using the above order on sets of integers).  

If an annulus
$A$ has minimal disks abutting it on both sides, we say that $A$ is
\underline{ambivalent}.
\end{defn}

\begin{lem} \label{lem:height} If $D_A$ is a minimal $\bdd$-compressing
disk for $A \in \Aaa \subset W$ then no component of $D_A \cap \Aaa$ is
an inessential arc in $\Aaa$.  In particular, $D_A$ is disjoint from all
spanning components of $\Aaa$.  Furthermore, if
$A'$ is another annulus in $\Aaa$ and $A'$ intersects $D_A$, then $A'$
is not ambivalent and any arc component of
$A' \cap D_A$ cuts off from $D_A$ a subdisk $D_A' \subset D_A$ that is
also a minimal disk.  
\end{lem}

\pf  Immediate, since otherwise replacing $D_A'$ by a minimal disk
(on the other side of $A'$ if possible) would decrease the complexity of
$D_A$.
\qed

\begin{defn} The side
opposite to the side on which a minimal disk abuts a non-spanning annulus
$A \in \Aaa$ is called the \underline{root side} of $A$.  If $A$ is
ambivalent, arbitrarily choose a side to call the root side.
\end{defn}  

\begin{defn} \label{defn:indep}  Let $\Aaa \subset \bdd_+ W$ be a
collection of essential non-spanning annuli in the boundary
of a compression body
$W$.  Then
$\Aaa$ is \underline{independently} \underline{longitudinal} if for each
$A
\in
\Aaa$ there is a meridian disk $D_A$ for
$W$ such that $\bdd D_A \cap \Aaa$ consists precisely of a single
spanning arc of $A$.
\end{defn}

An easy outermost arc argument shows that we can take the collection of
meridians that arise in Definition \ref{defn:indep} to be disjoint.

\begin{lem} \label{lem:glue}  Suppose $V$  and $W$ are  compression
bodies, and a collection of essential non-spanning annuli
$\Aaa \subset \bdd_+ W$ is independently longitudinal.  Then if $V$ is
attached to $W$ by identifying $\Aaa$ to any collection of annuli in
$\bdd_+ V$, the result is a compression body.
\end{lem}

\pf $\Aaa$ becomes a collection of boundary compressible annuli in
$V
\cup_{\Aaa} W$.  In this collection each annulus has a boundary
compressing disk of complexity $\{ 1 \}$ (i. e. a disk with interior
disjoint from
\Aaa).  Performing these boundary compressions turn $\Aaa$ into a
collection of disks.  So an alternate and equivalent construction for $V
\cup_{\Aaa} W$ would be to first compress $W$ along the meridians for
$\Aaa
\subset \bdd_+ W$ and then attach the resulting compression body $W'$ to
$\bdd_+ V$ along the disks in $\bdd_+ W'$ that are the remnants of
$\Aaa$. This alternate construction clearly gives a compression body. 
\qed

\begin{lem}  \label{lem:cbs} Let $P$ be a (not necessarily connected)
incompressible surface in a compression body $W$ such that $\partial P
\subset
\partial_+W$.  Then the result of cutting
$W$ along $P$ is a collection of compression bodies.
\end{lem}

\pf  See for instance \cite[Lemma 2]{Sc}. \qed

\begin{lem} \label{lem:push}
Suppose $\Aaa \subset W$ is a set of essential annuli in the compression
body $W$. Let
$\Aaa' \subset \Aaa$ be the subcollection of non-spanning annuli.  Cut out
from $W$ a collection of tunnels, one parallel to a spanning arc of each
annulus $A$ in $\Aaa'$ and lying on the root side of $A$.  Call the
resulting compression body $W'$.  Then the closure $W'_0$ of any
component of $W' -
\Aaa'$ is a compression body on which $\Aaa' \cap \bdd W_0$ is a
collection of independently longitudinal annuli.
\end{lem}

\pf  Let $A$ be any annulus in $\Aaa' \cap \bdd W'_0$ and let $W_0$
denote the component of $W - \Aaa$ from which tunnels were removed to get
$W'_0$.  If, in $W$, $W_0$ lay on the root side of $A$ then in
constructing $W'_0$ a tunnel has been removed that runs parallel to a
spanning arc of $A$.  The disk that defines the parallelism is a meridian
for $W'_0$ that intersects $A$ precisely in a spanning arc, as required. 

If, on the other hand, $W_0$ is not on the root side of $A$ then consider
a minimal $\bdd$-compressing disk
$D_A$ for $A$ in $W$.  The component $D_0$ of $D_A - \Aaa$ that abuts $A$
lies in $W_0$ and abuts any other annulus in $\Aaa$, if at all, only on
its root side (see Lemma \ref{lem:height}).  In particular, for each such
annulus, the tunnels that are removed to create $W'_0$ can be positioned
so that the component
$D'_0$ of $D_0 \cap W'_0$ that abuts $A$ runs over the tunnels instead of
across the other annuli.  See Figure 3.  $D'_0$ is then a meridian for
$W'_0$ that intersects
$\Aaa$ only in a single spanning arc of
$A$, as required.  \qed

\begin{figure}
\centering
\includegraphics[width=.8\textwidth]{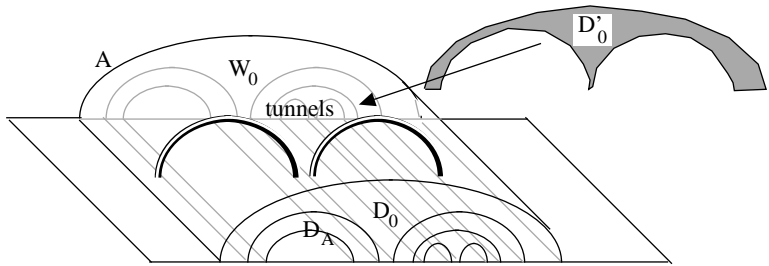}
\caption{Placing $D'_0$ as a meridian of $W'_0$.}
\end{figure} 

\section{Breaking up generalized Heegaard splittings by transverse annuli}
\label{sec:break}

A tunnel system for a knot corresponds to a Heegaard splitting of the
knot complement.  A Heegaard splitting, in turn, can be
untelescoped to produce a strongly irreducible generalized Heegaard
splitting, and vice versa.  When we consider degeneration of tunnel
number, we shall be concerned with constructing generalized Heegaard
splittings for each
$C(K^j), j = 1,...,n$ from a strongly irreducible generalized Heegaard
splitting for
$C(K^1 \# ... \# K^n)$.

The more general context for this section is this:   Let
$M$ be a compact orientable $3$-manifold and $\Aaa$ a properly imbedded
collection of essential annuli in $M$.   Suppose
$$(V_1 \cup_{S_1} W_1) \cup_{F_1} \dots \cup_{F_{m-1}} (V_m \cup_{S_m}
W_m)$$ is a strongly irreducible generalized Heegaard splitting of $M$
isotoped so that $(\Fff \cup \Sss) \cap A$ consists only of curves
essential in both $\Fff \cup \Sss$ and $\Aaa$, and such that this
number is minimal.

\begin{defn}  \label{defn:dipping} An annulus
component (or its closure) $A$ of $\Aaa - (\Fff \cup \Sss)$ is called a
\underline{dipping} annulus if, for some $1 \leq i \leq m$, $A \subset
V_i$ and
$\partial A \subset S_i = \partial_+V_i$.  
\end{defn}

\begin{thm} \label{thm:bif} Let $M^1,...,M^n$ be the components into
which $M$ is divided by the family of annuli $\Aaa$.
Then for each $1 \leq j \leq n$ there is a generalized Heegaard
splitting $(V_1^j \cup_{S_1^j} W_1^j)
\cup_{F_1^j} \dots \cup_{F_{m-1}^j} (V_m^j \cup_{S_m^j} W_m^j)$ of $M^j$
such that $$\sum_{i=1}^m J(V_i) \geq \sum_{j=1}^n \sum_{i=1}^m J(V_i^j)
 - 2k,$$ where $k$ is the number of dipping annuli among the components
of 
$\Aaa - (\Fff \cup \Sss)$.
\end{thm}

\pf The central problem in discerning Heegaard splittings in the
$M^j$ is that cutting
$V_i$ or
$W_i$ along
$\Aaa$ does not necessarily create compression bodies. For example,
cutting $\Fff \cup \Sss$ along $\Aaa$ does not even produce closed
surfaces.  We endeavor to remedy this fact by longitudinally attaching a
solid torus (which we could view as a collar of $\Aaa^j = \Aaa \cap M^j$)
to each of the annuli in $\bdd M^j$ and imbedding in each torus certain
annuli with longitudinal boundary.  These annuli, when attached to 
surfaces $S_i
\cap M^j$ (suitably stabilized) and $F_i \cap M^j$ will be shown to yield
generalized Heegaard splittings of the $M^j$.

The first step will be to describe how the new annuli are to be imbedded
in the solid torus collar $T^A$ of each $A \in \Aaa^j \subset \bdd M^j$. 
In $M$ itself, it's natural to define the "distance'' between two of the
surfaces in $(\Fff \cup \Sss) \subset M$ as the smallest number of
compression bodies one needs to pass through to get from a point in one
surface to a point in the other.  So, for example, the distance from
$\bdd_- V_1 = F_0$ to
$F_i$ is
$2i$ and from $F_0$ to $S_i$ the distance is $2i - 1$. 

$(\Fff
\cup \Sss) \cap A$ is a collection of parallel essential curves in the
annulus $A$. Let $\alpha$ be a spanning arc of $A$ that meets each
component of $(\Fff \cup \Sss) \cap A$ exactly once. Parameterize \aaa\ by
$0 \leq t \leq 1$. That is, choose a homeomorphism $h_1: \aaa \map [0,
1]$. Let $h_2: \aaa \map [0, 2m]$ be a continuous extension of the
function that assigns to each point in $\aaa \cap (\Fff \cup \Sss)$
its distance from $F_0$. We may as well take $h_2$ to be as simple as
possible. For example, on a segment of
$\aaa$ that runs between $F_i$ and $S_{i+1}$, say, define $h_2$ to
monotonically run from $2i$ to $2i + 1$.  On a segment of $\aaa$ that
lies in $W_i$ and has both ends on $S_i$,  define $h_2$ so that it has a
single maximum.

Let $h = (h_1, h_2): \aaa \map [0,1] \times [0, 2m]$ be the corresponding
imbedding of $\aaa$ in the first quadrant of
${\mathbb R}^2$, with the endpoints of $\aaa$ on the $x$-axis. Informally,
$h$ identifies \aaa\ with the graph of its distance from $F_0$.  Let
$D$ be the disk that lies below $h(\aaa)$ in $[0,1] \times [0, 2m]
\subset {\mathbb R}^2$.

Each line $y = 2i$ (respectively $y = 2i-1$) intersects $h(\aaa) \subset
{\mathbb R}^2$ at points at which $F_i$ (respectively $S_i$) intersects
\aaa.  In view of this correspondence, let ${\mathcal L}_{F_i}$ denote the
intersection of the line $y = 2i$ with $D$,  ${\mathcal L}_{S_i}$ denote the
intersection of the line $y = 2i-1$ with $D$, ${\mathcal L}_F = \cup_i {\mathcal
L}_{F_i}$ and ${\mathcal L}_S = \cup_i {\mathcal L}_{S_i}$.  See Figure 4.
Consider a component $R$ of $D - ({\mathcal L}_F \cup
{\mathcal L}_S)$.  $R$ is a polygon with an even number of sides.  The sides
lie alternately in $h(\aaa)$ and
${\mathcal L}_F \cup {\mathcal L}_S$.  

\begin{figure}
\centering
\includegraphics[width=.8\textwidth]{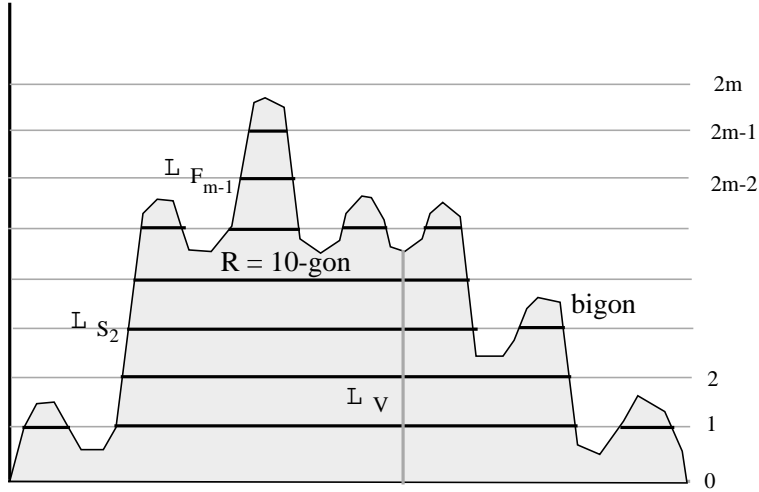}
\caption{The lines ${\mathcal L}_F, {\mathcal L}_S, {\mathcal L}_V
\subset V$.}
\end{figure} 

If $R$ is a bigon, one of its sides is a subarc of $h(\aaa)$ and one
of its sides is a component of ${\mathcal L}_S$, since an incompressible
annulus with both ends on $\bdd_- H$ in a compression body $H$ is
$\bdd$-parallel. Moreover $h_2$ has a maximum on the corresponding subarc
of $\aaa$, since
$D$ lies below the graph. That is, the annulus on which the corresponding
arc of \aaa\ lies is in some $W_i$.

If $R$ is a quadrilateral, two opposite sides are subarcs of
$h(\aaa)$, one side is a component of ${\mathcal L}_F$ and one a
component of
${\mathcal L}_S$.  If $R$ has $n$ sides with $n > 4$, then
$\frac{n}{2}$ sides will be subarcs of $h(\aaa)$, one side will be a
component of ${\mathcal L}_F$ and all other sides will be components of
${\mathcal L}_S$.  Thus there will be $\frac{n}{2} - 1$ sides that are
subarcs of
$h(\aaa)$ and that connect two sides that are components of ${\mathcal
L}_S$; these correspond to spanning arcs of dipping annuli in some $V_i$. 
Let $\Gamma$ be the collection of subarcs of
$\alpha$ that, in the boundary of some $R$, connect two sides that
are components of ${\mathcal L}_S$.  Then there is a correspondence between
the components of
$\Gamma$ and the collection of dipping annuli that lie in $A$.  

Now set $T^A = D \times S^1$.  $T^A$ will be attached to each of the two
copies of $A$ in $\cup_j (\bdd M^j)$ by the obvious identification
$h^{-1} \times {\bf 1}_{S^1}$ of $h(\aaa) \times S^1$ with $A
\cong (\aaa \times S^1)$.  The surfaces $(\Fff \cup \Sss) \cap M^j$ can
then be extended into $T^A$ by just attaching the collection of annuli
$({\mathcal L}_F \cup {\mathcal L}_S) \times S^1$. This operation will not
necessarily create compression bodies, but it will do so if we first
modify
$\Sss$ as described below.

Stabilize $\Sss$ by attaching tubes, one running parallel to a spanning
arc on each dipping annulus and lying on the root side of the annulus
in the compression body $V_i$ in which it lies. Denote the resulting
generalized Heegaard splitting of
$M$, now of genus
$k$ higher than originally, by $$(V'_1
\cup_{S'_1} W'_1) \cup_{F_1} \dots \cup_{F_{m-1}} (V'_m \cup_{S'_m}
W'_m).$$  Denote $\cup_i S'_i$ by $\Sss'$.

It follows from Lemma \ref{lem:push} that each component of the
complement in $V'_i$ of the dipping annuli is a compression body on whose
boundary the collection of incident dipping annuli are independently
longitudinal.  
  
\bigskip

{\bf Claim:}  Attaching each $T^A, A \subset \bdd M^j$ to $M^j$ and
capping off  the surfaces $(\Fff \cup \Sss') \cap M^{j}$ by the annuli
$({\mathcal L}_F \cup {\mathcal L}_S) \times S^1 \subset T^A$ yields a generalized
Heegaard splitting of $M^j$.

\bigskip

For any $i$, cutting along non-spanning components of $A \cap V'_i$ (or
$A \cap W'_i$) yields compression bodies by Lemma
\ref{lem:cbs}. Cutting along spanning annuli yields $(Q \times I) \cup
(1-handles)$ for some compact orientable surface $Q$.  When $T^A$ is
attached to $M^j$, then for $R$ a region as above, manifolds of the form
$R \times S^1$ are attached to $(Q \times I) \cup (1-handles)$.  

\vspace{1 mm}

\noindent {\bf Case 1:} $R$ is a bigon.

Note that a bigon corresponds precisely to a non-spanning annulus
component of $A \cap W_i$.  Then $R \times S^1$ is a solid torus that is
attached to $(Q \times I) \cup (1-handles)$ along a longitudinal
annulus of $R \times S^1$.  This does not change the homeomorphism type
of $(Q \times I) \cup (1-handles)$.

\vspace{1 mm}

\noindent {\bf Case 2:} $R$ is a quadrilateral. 

Attaching $R \times S^1$ to $(Q \times I) \cup (1-handles)$ yields $(Q'
\times I) \cup (1-handles)$, where
$Q'$ is the compact surface obtained by connecting two boundary components
of $Q$ by an annulus.

\vspace{1 mm}

\noindent {\bf Case 3:} $R$ is an $n$-gon with $n > 4$.  

Here attaching $R \times S^1$ to $(Q \times I)
\cup (1-handles)$ has the same effect as attaching a $(quadrilateral)
\times S^1$, but in addition, attachments are also made along dipping
annuli in the corresponding $V_i$.  $\Sss'$ has been constructed so that
each dipping annulus is independently longitudinal in the component of
$V'_i - \Aaa$ on which it lies, so the result is still of the form $(Q'
\times I) \cup (1-handles)$, essentially by Lemma
\ref{lem:glue}.

\vspace{1 mm}

Since all components of $\partial Q$ are eventually connected by annuli in
this process, the result is a union $V_i^j$ (or $W_i^j$) of compression
bodies.  \qed

\begin{rem} Theorem \ref{thm:bif} will eventually be applied to the family
of annuli $\Aaa = \{ A^1,...,A^{n-1} \}$ that decompose $C(K^1 \# \ldots
\# K^n)$ into $C(K^1), \ldots , C(K^n).$
\end{rem}

Certain annuli in the tori $T^A$, now described, will be useful in the
next section.

\begin{defn} \label{defn:plumb}
Let $D \subset {\mathbb R}^{2}$ be the disk described in the proof of Theorem
\ref{thm:bif}.  Let ${\mathcal L}_v$ be the intersection of a vertical line
$x = x_0$ with $D$.  Then the annulus ${\mathcal L}_v \times S^1$ is called
a \underline{plumbline} annulus in $T^A = D \times S^1$.  See Figure 4.
\end{defn}

\begin{rem} Let $M_+^j \cong M^j$ denote the manifold obtained from
$M^j$ by attaching the solid tori $\{ T^A| A \subset \bdd M^j \}$ to
$M^j$. Note that each plumbline annulus in $T^A \subset M_+^j$ will
intersect each splitting surface $S_i^j$ in at most one component.
\end{rem}

\bigskip

If we replace the collection of annuli in Theorem \ref{thm:bif} by a
collection of essential (= incompressible and not $\bdd$-parallel) tori,
then an analogous proof still applies.  The spanning arc $\alpha$ must be
replaced by an appropriate essential curve on $T$.  More specifically,
note that the components of
$(\Fff \cup
\Sss) \cap T$ are all parallel.  The curve that replaces $\alpha$ must be
an essential curve on $T$ that intersects each component of $(\Fff \cup
\Sss) \cap T$ exactly once.  The graph of this curve in $S^1 \times {\mathbb
R}^+$, constructed in analogy to $h(\aaa)$ in the proof of Theorem
\ref{thm:bif}, will cut out an annulus.  This annulus replaces $D$ in the
construction.  See Figure 5.  This yields the following result.

\begin{figure}
\centering
\includegraphics[width=.8\textwidth]{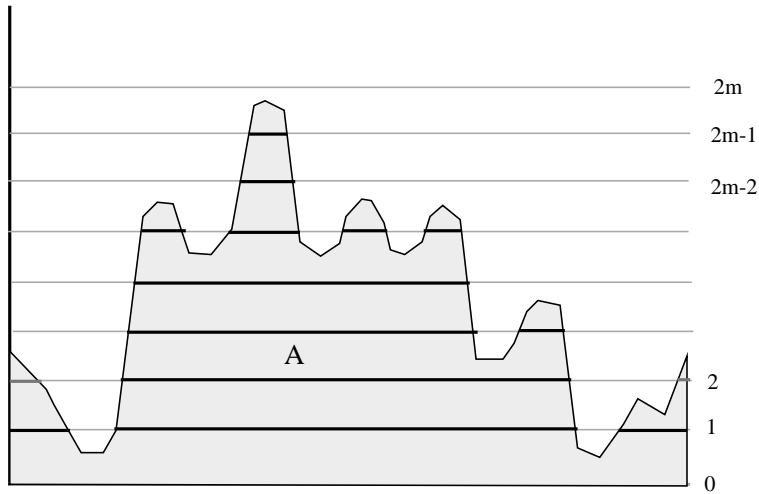}
\caption{Identify left and right sides so the $x$-axis becomes $S^1$.}
\end{figure}

\begin{thm} \label{thm:torus} Let $M$ be an orientable $3$-manifold
containing a family of essential tori $\mathcal
T$.  Suppose $(V_1 \cup_{S_1} W_1) \cup_{F_1} \dots \cup_{F_{m-1}} (V_m
\cup_{S_m} W_m)$ is a strongly irreducible generalized Heegaard splitting
of $M$ isotoped so that $(\Fff
\cup \Sss) \cap \mathcal T$ consists only of curves essential in both $\Fff
\cup \Sss$ and $\mathcal T$, and such that this number is minimal.  
Let $M^1,...,M^n$ be the components into which $M$ is divided by
$\mathcal T$. 

Then for each $1 \leq j \leq n$ there is a generalized Heegaard splitting
$(V_1^j
\cup_{S_1^j} W_1^j) \cup_{F_1^j} \dots \cup_{F_{m-1}^j} (V_m^j
\cup_{S_m^j} W_m^j)$ of $M^j$ such that $$\sum_{i=1}^m J(V_i) \geq
\sum_{j=1}^n \sum_{i=1}^m J(V_i^j)
 - 2k,$$ where $k$ is the number of dipping annuli in ${\mathcal T}
\cup (\cup_i V_i)$.
\end{thm}

\section{Destabilizations}
\label{sec:des}

\vspace{5 mm}

The generalized Heegaard splitting of each $M^j$ constructed in Theorem
\ref{thm:bif} may not be strongly irreducible and, even if it is, it
might still be simplified. Recall

\begin{defn} A Heegaard splitting $M = V \cup_S W$ is
\underline{stabilized} if there are properly imbedded disks $D_1 \subset
V$ and $D_2 \subset W$ such that
 $\vline \partial D_1 \cap \partial D_2 \vline = 1$.  
\end{defn}

\begin{rem} \label{rem:stab} In this case, cutting $V$ along $D_1$
(or $D_2$) yields another Heegaard splitting of lower genus.
\end{rem}

A strongly irreducible splitting may amalgamate to a stabilized or
reducible splitting, so (to account for this) the definition of a
stabilized generalized Heegaard splitting is necessarily a bit more
complicated.

\begin{defn} \label{defn:genstab}  A  generalized Heegaard splitting $M =
(V_1 \cup_{S_1} W_1)
\cup_{F_1} \dots \cup_{F_{m-1}} (V_m
\cup_{S_m} W_m)$ is \underline{stabilized} if there are disks $D_1$ and
$D_2$ such that $D_1$ is properly embedded in $(V_1 \cup_{S_1} W_1)
\cup_{F_1} \dots \cup_{F_{i-1}} V_i$, for some
$i$, and $D_2$ is properly embedded in $W_i \cup_{F_i} \dots
\cup_{F_{m-1}} (V_m \cup_{S_m} W_m)$.  Furthermore,
$\vline \partial D_1 \cap \partial D_2 \vline = 1$ and each component of
$D_j - (\Fff \cup \Sss)$ is either a disk or an annulus spanning the
compression body in which it lies.  
\end{defn}

If a generalized Heegaard splitting is stabilized (i. e. satisfies
Definition \ref{defn:genstab}) the associated amalgamated Heegaard
splitting is stabilized. It is easy to see that we may assume the annuli
components of
$D_j - (\Fff \cup \Sss)$ are essential.  It furthermore follows from
\cite{CG} that if the generalized splitting is strongly irreducible, then
we may assume each annulus component of $V_i \cap D_k$ or $W_i \cap D_k$
is a spanning annulus in a component of $V_i$ or $W_i$ that is a trivial
compression body.  

Note that if a generalized splitting is stabilized we can create a
generalized Heegaard splitting of lower genus (i. e. for which $\sum_i
J(V_i)$ is reduced) by amalgamating, reducing the genus as in Remark
\ref{rem:stab}, and then untelescoping again.

Here is an ad hoc criterion, useful in the present context, for showing
that a given generalized splitting is stabilized.  

\begin{defn}  \label{defn:tubed}

A surface $G \subset (D^2 \times S^1)$ is a \underline{tubed product} if
it can be described as follows:  Start with a properly imbedded
$1$-manifold ${\mathcal L} \subset D^2$, and a collection $\tau$ of $t$
disjoint arcs in the interior of $D$ such that $\tau \cap {\mathcal L} =
\bdd
\tau$.  Then
${\mathcal L} \times S^1 \subset D^2 \times S^1$ is a union of tori and
annuli.  Create $G$ by attaching unnested tubes to
${\mathcal L} \times S^1$ along the arcs $\tau \times \{ point \} \subset
D^2
\times
\{ point
\}$.  

Viewed dually, ${\mathcal L} \times S^1$ is obtained from $G$ by
simultaneously compressing $t$ disks.  For $G_0$ a component of $G$, each
component of ${\mathcal L} \times S^1$ that results from compressing $G_0$ is
said to {\em come from} $G_0$.

\end{defn}

\begin{lem} \label{lem:adhoc}

Suppose that a generalized Heegaard splitting  $(V_1 \cup_{S_1} W_1)
\cup_{F_1} \dots \cup_{F_{m-1}} (V_m \cup_{S_m} W_m)$ of a compact
$3$-manifold $M$ intersects an essential solid torus $D^2 \times S^1
\subset interior(M)$ so that the surface $G = (\Fff \cup \Sss) \cap (D^2
\times S^1)$ can be described as a tubed product with $t$ tubes as in
Definition \ref{defn:tubed}.  Suppose further that, for each $1 \leq i
\leq m$, at most one arc component of
${\mathcal L}
\times S^1$ comes from $S_i \cap (D^2 \times S^1)$.  Then the splitting is
stabilized at least $t$ times.
\end{lem}

\pf The proof is by induction on $|{\mathcal L}|$.  If
${\mathcal L} = \emptyset$ then there is nothing on which to attach
tubes, so $t = 0$ and there is nothing to prove. We may as well also
assume that no component of $\tau$ is parallel to a subarc of ${\mathcal
L}$ (that is, no disk component of $D - ({\mathcal L} \cup \tau)$
has boundary the union of a component of $\tau$ and a subarc of ${\mathcal
L}$).  For the corresponding tube is clearly a stabilization, so it
can be removed without affecting the truth of the lemma.

Suppose
$\mathcal L$ consists entirely of arcs.  Since each $F_i$ is
incompressible in both
$W_{i}$ and $V_{i+1}$, tubes can only have been attached to  components of
${\mathcal L} \times S^1$ coming from the thick surfaces $\Sss$.  By
assumption, there is at most one component of ${\mathcal L} \times S^1$ 
coming from any given $S_i$ and distinct $S_i$'s are separated by
components of
$\mathcal F$.  It follows that at least one component of $\tau$ is
parallel to a subsegment of ${\mathcal L}$, a contradiction.  See Figure
6.

\begin{figure}
\centering
\includegraphics[width=.5\textwidth]{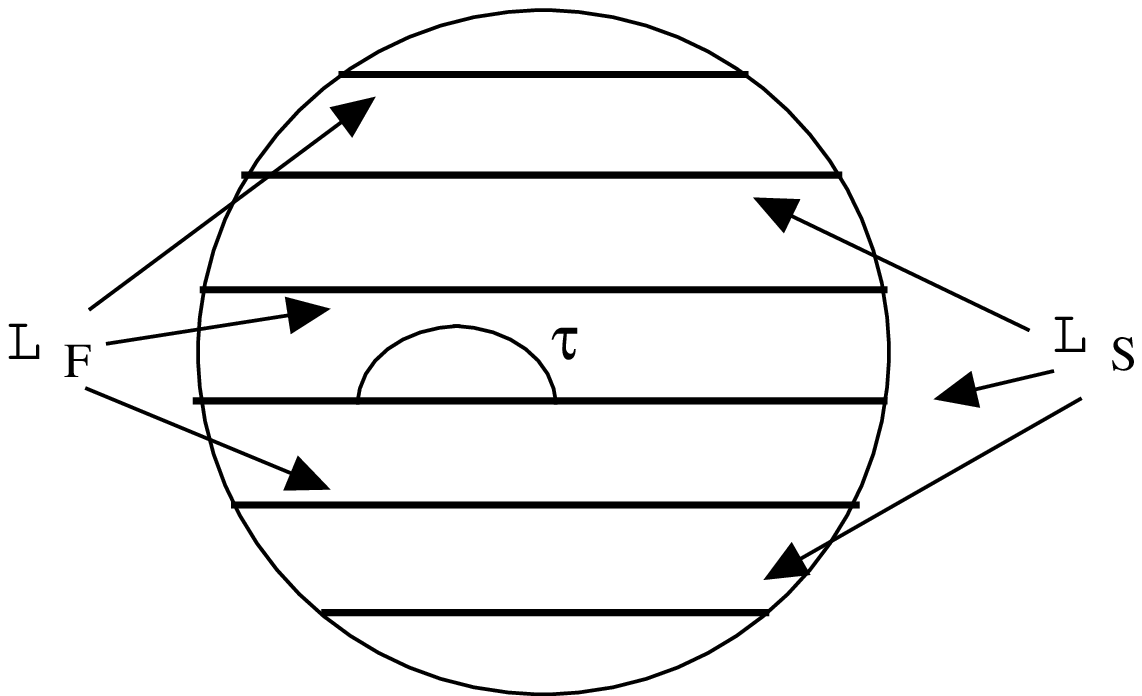}
\caption{}
\end{figure}

We proceed next to the case where there are closed components of $\mathcal
L$.  Each closed component $\llll \subset {\mathcal L} \subset D^2$ abuts two
components of $D^2 - \mathcal L$; call them the {\em inside} and {\em
outside} regions neighboring $\llll$, depending on whether or not $\llll$
separates the region from $\bdd D^2$.  If the outside neighboring region
of each closed component is an annulus, and the annulus spans some
$V_i$ or
$W_i$ then it follows, much as when there are no closed components, that
each tube attached is a stabilization.  So we can focus on a closed
component
$\llll \subset \mathcal L$ which is innermost among those whose outside
neighboring component is not a spanning annulus.  It follows that the
components of $D^2 - \mathcal L$ lying inside \llll\ consist exactly of
spanning annuli in the compression bodies, together with a single
compressing disk for some $S_i$. Then all tubes coming from arcs in
$\tau$ lying inside \llll\ in $D^2$ are stabilizations.  This means that
all components of $\Fff \cup \Sss$ lying inside the solid torus $U$
bounded by $\llll
\times S^1$ can be removed, and both $t$ and the number of stablizations
will be reduced by the number of tubes lying therein.  Indeed, $\llll
\times S^1$ itself can be removed, if it comes from $\Fff$.  

So, by induction, we are left with the case in which $\llll$ comes from
some $S_i$, cuts off from $D^2$ a compressing disk in $V_i$ or $W_i$, say,
$V_i$, and the outside neighbor of $\llll$ is not a spanning annulus of
$W_i$.  By judicious choice of curves satisfying these properties we may
further assume that no arc in $\tau$ has both ends on
\llll.   (For such an arc is either parallel to a subarc of \llll, and so
violates our assumption above, or cuts off from the outside neighboring
region another family of circles among which a substitute $\llll$ can be
found.  See Figure 7.)  If no tube is attached to \llll\ in the outside
neighboring region, the component $P$ of $W_i$ containing that region
will be a product, and so
$\bdd P \cup \bdd U$ can be removed from $\Fff \cup \Sss$ and still leave
a generalized Heegaard splitting.  

\begin{figure}
\centering
\includegraphics[width=.5\textwidth]{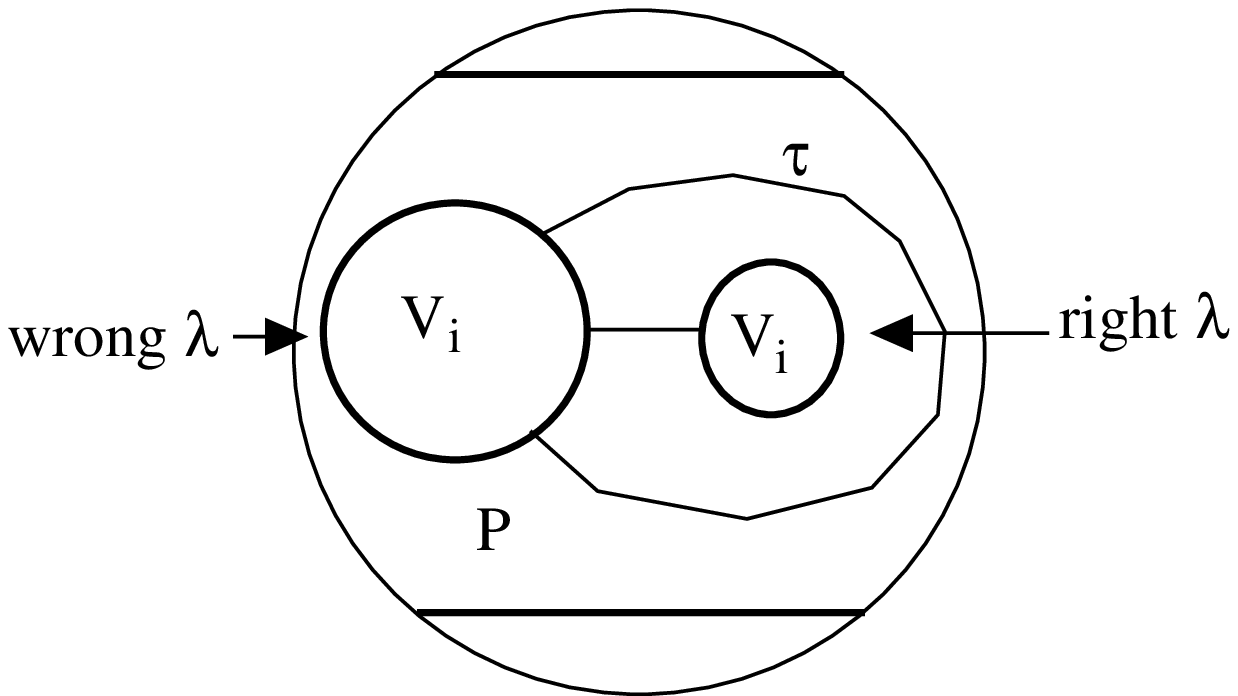}
\caption{}
\end{figure}

So we may as well assume that a tube is attached along an  arc
\aaa\ in the outside neighboring region of
\llll, and one end of \aaa\ lies on $\llll$ and the other end on another
component $\llll'$ of $\mathcal L$.  Necessarily $\llll'$ comes from $S_i$,
but possibly $\llll'$ is an arc component of $\mathcal L$. Whether $\llll'$
is an arc or is closed, the complement  of the tube in the annulus
$\aaa \times S^1$ is a disk which, together with a meridian disk
$\mu$ of $U$ form a stabilizing pair for the splitting.  Compressing
$S_i$ along $\mu$ then leaves a Heegaard splitting, still intersecting
$\Fff \cup \Sss$ in a tubed product, but with both $t$ and the number of
stabilizations reduced by one.  Since also $|{\mathcal L}|$ is reduced by
one, the result follows by induction.
\qed

\bigskip

We now turn to the setting of Theorem \ref{thm:bif} in which $M$ has a
strongly irreducible generalized Heegaard splitting  and will henceforth
assume that the family of annuli
$\Aaa$, separating
$M$ into components
$M^1,...,M^n$, is {\em complete}.  That is, if $(A', \bdd A') \subset (M,
\bdd M)$ is a properly imbedded incompressible annulus disjoint from
\Aaa, then $A'$ is $\bdd$-parallel in the component $M^j$ in which it
lies.  For example, if $M = C(K^1 \# ... \# K^n)$, then the knot
summation defines a collection $\Aaa$ of $n-1$ essential annuli in $M$
and this is a complete collection if and only if each knot is prime. 
Indeed, an annulus $A'$ in a knot complement $C(K^j)$, with $\bdd A'$ a
pair of meridian disks, can be extended to become a decomposing sphere in
$S^3$ by capping off $\bdd A'$ with meridian disks.  A resulting summand
is trivial if and only if $A'$ is  $\bdd$-parallel in $C(K^j)$.

Note that the condition that $\Aaa$ is complete is weaker than the
assumption that each $M^j$ is acylindrical, since it says nothing about
incompressible annuli in $M^j$ whose boundaries cross the curves $\bdd
\Aaa \cap \bdd M^j$. 

\begin{defn} \label{defn:wide}  Suppose that the surface $\Fff \cup \Sss$
has been isotoped to intersect $\Aaa$ only in curves essential in both
$\Fff
\cup \Sss$ and $\Aaa$ and so  that the number of components of $(\Fff
\cup \Sss) \cap \Aaa$ is minimal subject to this condition.  Then an
annular component $\ooo$ of $\Sss - \Aaa$ is called a \underline{wide}
annulus if it is adjacent to (that is, shares a boundary component with)
a dipping annulus in $\Aaa - (\Fff \cup \Sss)$.  

A component $U$ of $W_i - \Aaa$ (resp $V_i - \Aaa$) is 
\underline{exceptional} if, for an annulus $B$, there is an embedding
$(B,
\bdd B) \times I
\subset W_i$ and an open disk $D$ in $B$ so that $U$ is the image
of $(B - D) \times I$, with $U \cap \Aaa = \bdd(B) \times
I$.  
\end{defn}

Exceptional components were first described near the end
of the proof of Theorem  \ref{thm:indic}. See Figure 1.  Note that $U
\cap S_i$ is a
$4$-punctured sphere and the closure of $D
\times \{ point \}$ is a compressing disk for $U \cap S_i$ in
$V_i$ (resp $W_i$).  

\begin{defn} For $U$ an exceptional component of $W_i - \Aaa$ (resp $V_i
- \Aaa$), the $4$-punctured sphere $U_S = \bdd U \cap S_i$ will be
called an \underline{exceptional} component of $S_i - \Aaa$ and the
annuli obtained by compressing $U_s$ into $V_i$ (resp $W_i$) will be
called \underline{virtual annuli} of $S_i$.  
\end{defn}

\begin{lem} \label{lem:dest} Let $M$ be a compact
$3$-manifold with strongly irreducible generalized Heegaard 
splitting $M = (V_1 \cup_{S_1} W_1) \cup_{F_1}
(V_2 \cup_{S_2} W_2) \cup_{F_2} \dots \cup_{F_{m-1}} (V_m \cup_{S_m}
W_m)$.  Let $M^1,...,M^n$ be the components into which $M$ is divided by
the complete collection of annuli $\Aaa$.  Suppose there
are
$w$ wide annuli among the components of $\Sss - \Aaa$ and $e$ exceptional
components.  Suppose further that each Heegaard
splitting
$(V_1^j
\cup_{S_1^j} W_1^j)
\cup_{F_1^j} \dots \cup_{F_{m-1}^j} (V_m^j \cup_{S_m^j} W_m^j)$ of $M^j$
in the construction in Theorem \ref{thm:bif} can be destabilized
$d_j$ times.  Then $\Ss_j d_j \geq \frac{w}{2} + e$.
\end{lem}

\pf  Let $B \subset M^j$ be an annulus disjoint from $\Fff \cup
\Sss$ whose ends are essential (hence core curves) in components $A_0$
and $A_1$ of \Aaa.  Here we will allow $A_0 = A_1$ but not if $B$ is
parallel to a subannulus of $A_0 = A_1$. Unless $\Sss \cap M^j$ contains
no  annuli or exceptional components, such an annulus can be found, e. g. 
parallel to an annulus component of $\Sss \cap M^j$ (or a virtual
annulus).  Since $\Aaa$ is complete, $B$ is boundary parallel in $M^j$ and
the annulus $B^+_{\bdd} \subset \bdd M^j$
to which it is parallel contains as collars of its ends subannuli
$A'_0
\subset A_0$ and $A'_1 \subset A_1$ with $A'_0 \neq A'_1$. Denote the
annulus $B^+_{\bdd} \cap \bdd M = B^+_{\bdd} - (A'_0 \cup A'_1)$ by 
$B_{\bdd} \subset \bdd M$. It is known (see for example \cite{Sh}) how
a strongly irreducible Heegaard splitting can intersect a solid torus
under the conditions here, so we know that each $S_i$
intersects the solid torus $T$ lying between $B$ and $B^+_{\bdd}$ in
$M^j$ in a collection of $\bdd$-parallel annuli, plus possibly a
component in which two such annuli are tubed together by a
$\bdd$-parallel tube, i. e. an exceptional component.

Since the number of curves in $\Aaa \cap (\Fff
\cup \Sss)$ has been minimized, none of the annuli
of $(\Fff \cup \Sss) \cap T$ (and none of the virtual annuli) has both its
ends on the same component of $\bdd T \cap \Aaa$.  Thus each annulus
component of $(\Fff \cup \Sss) \cap T$  (or virtual annulus) is an annulus
much like $B$. So with no loss of generality, we may as well assume that
$B$ is outermost, i. e. no other annulus component of
$(\Fff \cup \Sss) \cap M^j$ (or virtual annulus) cuts off a solid torus
containing $B$.  Also, $T$ can be parameterized as $disk \times S^1$ so
that the annuli of $(\Fff \cup \Sss) \cap T$ plus the virtual
annuli are just the product of a collection of proper arcs in $disk$ with
$S^1$.

Now attach to $M^j$ the solid tori described in Theorem \ref{thm:bif} to
get $M_+^j$. Assume, for initial simplicity, that neither end of $B$ abuts
a dipping annulus on its root side.  Extend $B$ by attaching to the ends
of
$B$ at $A_0$ and $A_1$ the plumbline annuli in the tori $T^{A_0}$ and
$T^{A_1}$ described in Definition \ref{defn:plumb}.  Call the resulting
annulus
$B_+$ and denote by $T_+$ the solid torus that $B_+$
cuts off from $M_+^j$.  By construction, the Heegaard splitting $(V_1^j
\cup_{S_1^j} W_1^j) \cup_{F_1^j}
\dots \cup_{F_{m-1}^j} (V_m^j \cup_{S_m^j} W_m^j)$ intersects $T_+$ in a
tubed product, where the tubes are attached in the manner described in
Theorem \ref{thm:bif} or, in the case of exceptional components, by the
compressions that create the virtual annuli. By Lemma
\ref{lem:adhoc} it suffices to show that the number of tubes  not coming
from exceptional components in this tubed product is at least half as
large as the number of wide annuli among the components of $\Sss \cap T$.
The argument is little changed by assuming, as we henceforth do, that
there are no exceptional components.

Under the homeomorphism $T \cong disk \times S^1$, any component $T_0$ of
$T - (\Fff \cup \Sss)$ is the product of a $2p$-gon with
$S^1$.  The sides of the $2p$-gon become annuli components lying in
alternately $(\Fff \cup \Sss \cup \bdd M) - \Aaa$ and in $\Aaa - (\Fff
\cup \Sss)$.  At most one of the sides in each $2p$-gon lies in $\Fff
\cup \bdd M$.  If $T_0$ lies in some $V_i$ then either $p-2$ or $p$ of
its sides are dipping annuli, depending on whether or not $T_0$ abuts
$\Fff \cup \bdd M$.  When $T_0$ abuts $\Fff \cup \bdd
M$ then the number of wide annuli is correspondingly $0$ if $p = 2$ and
$p-1$ if $p > 2$.  When $T_0$ does not abut $\Fff \cup \bdd M$ it is
simply $p$.  

It similarly follows from the
construction in Theorem \ref{thm:bif} that the number of tubes added in
$T_0$ to create $S^j$ is $p-2$ if $T_0$ abuts
$\Fff \cup \bdd M$ and is $p$ or $p-1$ if it does not.  The
required inequality then follows from the inequalities 
$p-2 \geq \frac{p-1}{2}$ when $p > 2$ and $T_0$ abuts $\Fff \cup \bdd M$
and $p-1 \geq \frac{p}{2}$ when $p > 1$ and $T_0$ does not abut $\Fff \cup
\bdd M$.

The previous argument can be extended to include the case in which one
or both ends of $B$ abut a dipping annulus $A$ on its root side (when the
exact construction would need to handle a tube running through $B$; cf.
Figure 8): Just slide the end of $B$ out (away from $T$) just beyond the
end of
$A$ (so the end of $B$ cuts off a small collar of a boundary component of
$\Sss
\cap \Aaa$) before attaching plumbline annuli in $T_A$. The argument then
goes through as above.  \qed

\begin{figure}
\centering
\includegraphics[width=0.8\textwidth]{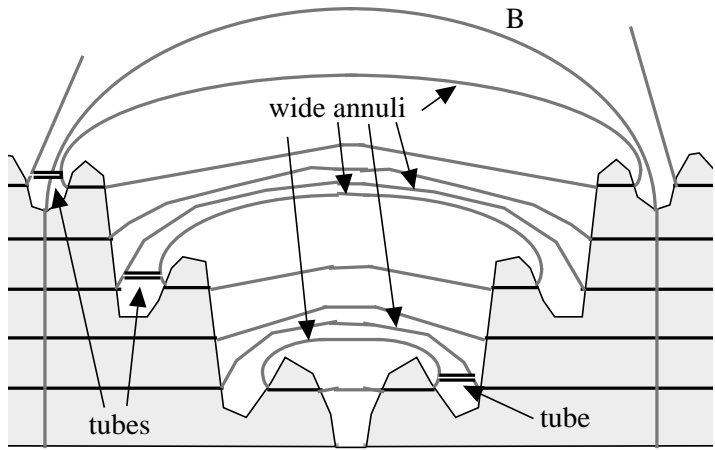}
\caption{Wide annuli and stabilizing tubes}
\end{figure} 

\section{Counting Annuli in the general case}
\label{sec:generalcount}

\vspace{5 mm}

In this section we use cut and paste techniques to bound
the number of dipping annuli in a strongly irreducible generalized
Heegaard splitting.  

\begin{defn} \label{def:parallel} Suppose $\Aaa \subset M$ is a properly
embedded collection of annuli.  Then a component of $M - \Aaa$ that is a
solid torus whose boundary intersects $\Aaa$ in a single longitudinal
annulus is a \underline{parallelism}.

\end{defn}  

Put another way, an annulus that is $\bdd$-parallel cuts off a
parallelism, unless there are other components of $\Aaa$ lying between it
and $\bdd M$.

\begin{lem} \label{lem:disk}   Suppose that $(\Aaa, \bdd \Aaa) \subset
(W, \bdd_+ W)$ is a collection of essential annuli in a compression
body $W$ and suppose $J(W) > 0$ (i. e. $W$ is not a product).  Then
there is an essential disk $D$ in $W$ such that $D$ is disjoint from
$\Aaa$.  Moreover, if not all components of $\bdd A$ are spanning, then
$D$ may be chosen so that at least one side of
$D$ in $W - \eta(D)$ lies on a parallelism of $\Aaa$ in $W - \eta(D)$. 
\end{lem} 

\pf  Let ${\mathcal D} \neq \emptyset$ be a set of defining disks for
$W$.  We argue by induction on $\vline \Aaa \cap {\mathcal D} \vline$. If
$\Aaa \cap \Daa = \emptyset$ then all non-spanning components of $\Aaa$
are $\bdd$-parallel in the product $W - \eta({\mathcal D})$.  One of them
cuts off a parallelism $U$ which, because $\Aaa$ is essential in $W$, is
adjacent to some of the disks in $\mathcal D$.  Then a disk in the boundary
annulus $\bdd U - \Aaa$ that contains all such disks corresponds, back in
$W$, to the required disk $D$.

If $\Aaa \cap \Daa \neq \emptyset$, consider ${\mathcal D} \cap \mathcal
A$.  Since the annuli are incompressible, we can easily remove any circle
of intersection.  If there is an arc  in
${\mathcal D} \cap \mathcal A$ that is inessential in $\mathcal A$, then let $\alpha$
be an outermost such arc in $\mathcal A$, and we may cut the disk $D$ in
$\mathcal D$ containing $\alpha$ along $\alpha$ and paste on two copies of the
disk cut off by $\alpha$ in $\mathcal A$ to obtain a new disk $D'$.  Replacing
$D$ by $D'$ in $\mathcal D$ produces a new set of defining disks $\mathcal D'$
with $\vline {\mathcal A} \cap {\mathcal D'} \vline < \vline {\mathcal A} \cap {\mathcal
D} \vline$.  The result follows by induction.

On the other hand, if all arcs in ${\mathcal A} \cap {\mathcal D}$ are essential
in $\mathcal A$, let $\beta$ be an arc in ${\mathcal A} \cap {\mathcal D}$ that is
outermost in $\mathcal D$.  Let $A$ be the annulus in $\mathcal A$ that gives
rise to
$\beta$. Cutting and pasting $A$ along $\beta$ and the outermost disk cut
off in $\mathcal D$ yields a disk $D'$ disjoint from $\mathcal A$.  Since $A$ is
essential, it follows that $D'$ is essential, and is as required. 
\qed

\begin{lem} \label{lem:generalcount}  Suppose there is a collection
$\Baa$ of essential annuli in a compression body $W$ and $\Aaa \subset
\Baa$ is the subcollection of non-spanning annuli. Let $W_- \subset W -
\eta (\Baa)$ consist of those components which are incident to $\bdd_-
W$.   Let $S = \bdd_+ W - \eta(\Baa)$, let
$S_- \subset S = S \cap W_-$ and let $S_+ = S - S_-$. Let $a$ denote the
number of annulus components of $S_+$.

Then $\Aaa$ has at most $J(W) + \frac{a}{2}$ components.
\end{lem}

\pf With no loss of generality we can ignore spanning annuli and
assume $\Baa = \Aaa$. 

Let $c(W, A) = 2J(W) + a - 2|\Aaa|$; it suffices to show that $c(W,A) \geq
0$.  The proof is by induction on $J(W)$. When $J(W) = 0$ there are no
annuli, and there is nothing to prove.   When $J(W) > 0$ it follows from
Lemma
\ref{lem:disk} that there is a $\bdd$-reducing disk $D$ for
$W$ that is disjoint from $\Aaa$.  The result of cutting
$W$ along $D$ is either one or two compression bodies $W'$ with $J(W')
=J(W) - 2$.  In particular, $J$ is lower in (each component of) $W'$, so
we can assume the Lemma is true in $W'$.
Remove all inessential annuli from
$\Aaa \cap W'$ and call the result $\Aaa'$.  Define $c' = c(W', \Aaa')$. 
By induction, it suffices to show that $c' \leq c = c(W,A)$. 

The first step, cutting $W$ along $D$, decreases $2J(W)$ by
$4$ and raises $a$ by the number of inessential annuli
that result, no more than $2$.  So this step decreases $c$ by at least
two. It may create, however, one or two components of $W' - \Aaa'$ that
are parallelisms.

Next examine what happens when an inessential annulus is removed:
$-2|\Aaa|$ goes up by two.  Also $a$ drops by at least one and will drop
by two if $A$ is adjacent to an annulus component of $\bdd_+ W - \Aaa$.
The latter will happen, for example, if the inessential annulus that's
removed is parallel to another one.  Thus this step either reduces the
number of components of $W' - \Aaa'$ that are parallelisms  and
simultaneously increases $c$ by at most one, or leaves both $c$ and the
number of parallelisms unchanged. Continue the process until all
parallelisms (which, at the beginning, are no more than two) are
eliminated.  The result is to increase $c$ by at most two.  Combining
both steps, $c$ has not increased.
\qed 

\bigskip

Return now to the original context:  $M$ is a compact orientable
$3$-manifold and $\Aaa$ a  complete properly imbedded
collection of essential annuli in $M$.   Suppose $$(V_1 \cup_{S_1} W_1)
\cup_{F_1} \dots \cup_{F_{m-1}} (V_m \cup_{S_m} W_m)$$ is a strongly
irreducible generalized Heegaard splitting of $M$ isotoped so that $(\Fff
\cup \Sss)
\cap A$ consists only of curves essential in both $\Fff
\cup \Sss$ and $\Aaa$, and such that this number is minimal.  Let
$M^1,...,M^n$ be the components into which $M$ is divided by the
family of annuli $\Aaa$. We apply Lemma
\ref{lem:generalcount} to the annuli $\Aaa \cap V_i$.  

\begin{thm} \label{thm:generalcount}
For $M$ and  $\Aaa$ as described above, let $k$ be the number of 
dipping annuli in $\Aaa - \Sss$, and $w$ be the number of wide annuli in
$\Sss - \Aaa$. Then
$$k \leq (\Ss_{i=1}^m J(V_i)) + \frac{w}{2}.$$
\end{thm}

\pf Lemma \ref{lem:generalcount} says that for each $1 \leq i \leq
m$ the number of dipping annuli in $\Aaa
\cap V_i$ is at most $J(V_i) + \frac{a}{2}$, where $a$ (defined
in Lemma \ref{lem:generalcount}) counts
the number of annuli components of $S_i - \Aaa$ with a property that
assures that they are adjacent to dipping annuli and therefore wide. 
That is, when Lemma
\ref{lem:generalcount} is applied to $V_i$ the inequality remains true if
$a$ is replaced by the number of wide annuli in $S_i$.  Summing over all
$V_i$ we get
$k \leq (\Ss_{i=1}^m J(V_i)) + \frac{w}{2}$ as required.  
\qed
 
In the appendix, Andrew Casson constructs an example that the inequality
of Theorem  \ref{thm:generalcount} is in some sense best possible. That
is, given any $c < 1$ there is an example for which $w = 0$ and $k >
c(\Ss_{i=1}^m J(V_i)).$

\begin{cor} \label{cor:general}
Let $M$ and $\Aaa$ be as described above.  Then there is a Heegaard
splitting $V_-^j \cup_{S_-^j} W_-^j$ for each $M^j$ so that
$$3(\Ss_{i=1}^m J(V_i)) \geq \sum_j J(V_-^j).$$
\end{cor}

\pf  Let $k$ denote the number of dipping annuli in $\Aaa \cap
(\cup_i V_i)$, and let $w$ denote the number of wide annuli in $\Sss -
\Aaa$.  The construction in Theorem \ref{thm:bif} yields generalized
Heegaard splittings for $C(K^j)$ such that $2k + \sum J(V_i) \geq
\sum_{i,j} J(V_i^j)$.  

Now substitute for $k$ from Theorem \ref{thm:generalcount} to get
$$3(\Ss_{i=1}^m J(V_i)) + w \geq \sum_{i,j}
J(V_i^j).$$  According to Lemma \ref{lem:dest} the induced 
Heegaard splittings for the $C(K^j)$ can be destabilized at
least $\frac{w}{2}$ times, yielding Heegaard splittings $V_-^j
\cup_{S_-^j} W_-^j$ for each $C(K^j)$ with $\sum_j J(V_-^j) \leq \sum_{i,j}
J(V_i^j) - w$.  It follows that 
$$3(\Ss_{i=1}^m J(V_i)) \geq \sum_j J(V_-^j).$$  \qed

\begin{cor} \label{cor:generaltunnel}
If $K^1, \ldots, K^n \subset S^3$ are prime knots
knots then $$t(K^1 \# \ldots \# K^n) \geq \frac{1}{3}(t(K^1) + \ldots +
t(K^n)).$$
\end{cor}

\pf  By Lemma \ref{lem:ind} there is a generalized Heegaard
splitting $$C(K^1 \# \ldots \# K^n) \cong (V_1 \cup_{S_1} W_1) \cup_{F_1}
\dots \cup_{F_{m-1}}  (V_m \cup_{S_m} W_m)$$ for which $\sum_{i=1}^m
J(V_i) = 2t(K^1 \# \ldots \# K^n)$. Let $\Aaa$ be the $n-1$ annuli that
decompose
$C(K^1 \# \ldots \# K^n)$ into the complements of the constituent prime
knots.   We may assume that $\Fff \cup \Sss$ and $\Aaa$ have
been chosen so that $(\Fff \cup \Sss) \cap \Aaa$ consists only of
essential curves and then also minimizes the number of such
intersections. The family $\Aaa$ is complete since each $K^j$ is prime. 
Apply Corollary \ref{cor:general}, substitute from Lemma
\ref{lem:ind}, and divide by two to get $3t(K^1 \# \ldots \# K^n)
\geq t(K^1) + \ldots + t(K^q)$ as required.  \qed

\bigskip
Following Theorem \ref{thm:torus}, the same sort of argument can be
applied to essential tori in $M$.  For $M$ a knot complement, the
application is to satellite knots.

\begin{defn}
\cite{BZ} Let $\tilde K$ be a knot in a $3$-sphere $S^3$ and $V$
an unknotted solid torus in $S^3$ with $\tilde K \subset V \subset
S^3$.  Assume that $\tilde K$ is not contained in a
$3$-ball of $V$.  A homeomorphism $h: V \rightarrow \hat V$
onto a tubular neighborhood $\hat V$ of a nontrivial knot
$\hat K \subset S^3$ which maps a meridian of $S^3 - V$ onto a
longitude of $\hat K$ maps $\tilde K$ onto a knot $K = h(\tilde K)
\subset S^3$.  The knot $K$ is called a
\underline{satellite} of $\hat K$, and $\hat K$ is called its
\underline{companion}.  The pair $(V, \tilde K)$ is called a
\underline{pattern} of $K$.
\end{defn}

\begin{thm} \label{thm:sat} Let $K$ be a satellite knot, then $t(K) \geq
\frac{1}{3}(t(\tilde K) + t(\hat K))$.
\end{thm}

\pf  We merely sketch the proof.  Let $(V_1 \cup_{S_1} W_1)
\cup_{F_1} \dots \cup_{F_{m-1}}  (V_m \cup_{S_m} W_m)$ be a generalized
Heeggard splitting for $C(K)$ for which $\sum_{i=1}^m J(V_i) = 2t(K)$.
Isotope $\bdd V$ so it intersects $\Fff \cup \Sss$ only in curves
essential in both, and in a minimal number of these curves.  Attach
copies of $torus \times I$ to $S^3 - \hat V$ and to $\hat V$ and complete
the surfaces $(\Fff \cup \Sss) - \hat V$ and $(\Fff \cup \Sss) \cap \hat
V$ to give Heegaard splittings of $C(\hat K)$ and $V - \eta(\tilde K)$
respectively, following Theorem \ref{thm:torus}. The latter can be made a
Heegaard splitting of $C(\tilde K)$ by just filling in a solid torus along
$\bdd V$.  Destabilizations can be found for these Heeegaard splittings
just as in Lemma \ref{lem:dest}.  \qed

\section{Counting annuli - specialized case}
\label{sec:specialcount}

We can improve the count of dipping annuli in Section
\ref{sec:generalcount} if we add the one further assumption that each
component of $\Sss - \Aaa$ has an even number of boundary components. 
This condition will certainly be satisfied when the annuli are the
decomposing annuli from a knot summation.  More generally, 

\begin{defn} \label{def:parityclass} For ${\mathcal C}$ a collection of
circles, the \underline{parity class} $\aaa \in H^1({\mathcal C}, {\mathbb
Z}_2)$ is the class that is non-trivial on each fundamental class of a
component of
$\mathcal C$.
\end{defn}

Informally, the parity class just counts the parity of the number of
components, even or odd.

\begin{lem} \label{lem:parity}  Suppose the properly embedded essential
annuli $\Aaa \subset M$ have the property that the parity class of
$\bdd
\Aaa \subset \bdd M$ is the restriction of a class in $H^1(M, {\mathbb
Z}_2)$. As usual, let $M^1, \ldots , M^j$ be the closed
complementary components of $\Aaa$ in $M$, so each annulus in $\Aaa$
becomes two annuli in $\bdd (\cup_j M^j)$. If
$(S, \bdd S) \subset (M^j, \Aaa \cap M^j)$ is an orientable, properly
embedded surface in any $M^j$, then $|\bdd S|$ is even.
\end{lem}

\pf
We may as well assume $S$ is connected. Let $[S]$ denote the
non-trivial class of $H_2(S, \bdd S; {\mathbb Z}_2)$, let \aaa\ be the
parity class of $\bdd \Aaa \subset \bdd M$ and let
$\tilde{\aaa}
\in H^1(M)$ be a class such that $i^*(\tilde{\aaa}) = \aaa$. Then,
after homotoping $\bdd S$ to $\bdd A$ we see that the evaluation
$[\aaa, \bdd [S]] = [\delta(\aaa), [S]] = [\delta(i^*(\tilde{\aaa})),
[S]] = 0$ by the exactness of the cohomology sequence for the pair $(M,
\bdd A)$.   \qed

\begin{defn} \label{def:paritycond}  A collection of annuli $\Aaa \subset
M$ as in Lemma \ref{lem:parity} satisfies the
\underline{parity condition}.
\end{defn}

This section will be a repeat of Section \ref{sec:generalcount} in the
special case in which the annuli $\Aaa$ satisfy the parity condition in
$M$, so every component of $\Sss - \Aaa$ has an even number of boundary
components.  The more delicate analysis will require a new notion:

\begin{defn} \label{def:special} Suppose that
$(\Aaa, \bdd \Aaa) \subset (W, \bdd_+ W)$ is a collection of essential
annuli in a compression body
$W$.  A component $U$ of $W - \eta(\Aaa)$ is called \underline{special}
if there is a planar surface
$P$ with boundary components
$p_0, ..., p_r$ and a homeomorphism $(U, \bdd U \cap \Aaa)
\cong (P \times I, \cup_{i=1}^{r} p_i \times I)$. See Figure 9.

The \underline{index}
$I(U)$ of the special component $U$ is defined to be $3 - r$.  
\end{defn}

\begin{figure}
\centering
\includegraphics[width=.6\textwidth]{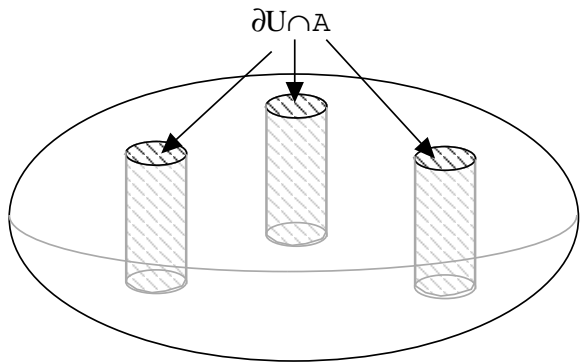}
\caption{A special component with $r = 3, I = 0$.}
\end{figure}

For example, a boundary parallel annulus cuts off from $W$ a
parallelism (see Definition \ref{def:parallel}).  This is a special
component of index
$2$, for in this case
$P$ is an annulus and $r = 1$. More generally, a special component is
obtained from a collection of
$r$ such examples by connecting them with $r-1$ $1$-handles. An
exceptional component (see Definition \ref{defn:wide}) is special in the
compression body in which it lies.  Its index is $1$. Any $\bdd$-reducing
disk $(D, \bdd D) \subset (U, \bdd U \cap \bdd_+ W)$ for $W$ that lies in
a special component $U$ divides $U$ into two special components in the
boundary-reduced manifold $W'$, and the sum of their indices is $index(U)
+ 3$. 

We now embark on improving the bound
of Lemma \ref{lem:generalcount} under the parity assumption: every
component of $\bdd W - \Baa$ has an even number of boundary components.
The argument is in three stages, incorporated in the following lemmas
that culminate in Lemma \ref{lem:supercount}:

\begin{lem} \label{lem:subcount}
Suppose there is a collection
$\Baa$ of essential annuli in a compression body $W$ and
$\Aaa \subset \Baa$ is the subcollection of non-spanning annuli.  Let $W_-
\subset W - \eta (\Baa)$ consist of those components which are incident
to $\bdd_- W$.   Let $S = \bdd_+ W - \eta(\Baa)$ and suppose every
component of $S$ has an even number of boundary components.  Let
$S_- \subset S = S \cap W_-$ and let $S_+ = S - S_-$. Let $m$
be the number of annuli in \Aaa\ that are adjacent to $W_-$ (counted
twice if $W_-$ abuts the annulus on both sides),
$n$ be the number of components of \Aaa\ that are not adjacent to $W_-$,
and
$p$ be the number of spanning annuli.  Let $s$ be the number of special
components of $W - \Baa$ (see Definition \ref{def:special}). 

Then  $$J(W) + |S_+| - |S_-| + genus (S_-) - m + p +
\frac{\chi(\bdd_- W)}{2} - 2n +  s \geq 0.$$
\end{lem}

\pf If a spanning component of $\Baa$ is removed, the only effect is to
lower
$p$ by one and to raise $[-|S_-| + genus (S_-)]$ by one.  So we may as
well assume that there are no spanning annuli, so $\Baa = \Aaa$ and $p =
0$.

The proof will be by induction on $J = J(W)$.  When $J = 0$ then $W$ is a
product, and there are no annuli.  In particular, $S_- \cong \bdd_- W$ so
the inequality follows from 
$$|\bdd_- W| - genus (\bdd_- W) = \frac{\chi(\bdd_- W)}{2}.$$

Adding $1$-handles on $\bdd_+ W$ only increases $J(W) + genus(S_-)$ and
has no other effect, so the Lemma follows from the case $J = 0$ whenever
$\Aaa = \emptyset$.  So we henceforth assume $\Aaa \neq \emptyset$.  

For the collection $(\Aaa, \bdd \Aaa) \subset (W, \bdd_+ W)$ define
$$b(W, \Aaa) = J(W) + |S_+| - |S_-| + genus (S_-) - m + p - 2n +  s.$$

Following Lemma \ref{lem:disk}, there is an essential disk $D$ in $W$
that is disjoint from $\Aaa$ and cuts off a parallelism from the
component of $W - \Aaa$ in which it lies.  The result of cutting
$W$ along $D$ is either one or two compression bodies $W'$ with $J(W')
=J(W) - 2$.  In particular, $J$ is lower in (each component of) $W'$, so
we can assume the Lemma is true in $W'$.  Recall that $D$ has been chosen
so that at least one of the adjacent components of $W' - \Aaa$ is a
parallelism of $\Aaa$ in $W'$.   Remove all inessential annuli from
$\Aaa \cap W'$ and call the result $\Aaa'$.  

Abbreviate $b(W, \Aaa)$ to simply $b$ and set $b' = b(W', \Aaa')$.  We
will show that $b \geq b'$. Since it follows from the inductive
hypothesis that $$b' \geq  -\frac{\chi(\bdd_- W')}{2} = 
-\frac{\chi(\bdd_- W)}{2},$$ we will then be able to conclude that 
$$b \geq  -\frac{\chi(\bdd_- W)}{2}$$ as required.

The process by which $b'$ can be calculated from $b$ consists of two
steps:  First $\bdd$-reduce $W$ along $D$, then remove all resulting
inessential annuli.  We examine the effect of each move on $b$ (in the
interim extending the definition of $b$ also to the case in which $\Aaa$
may contain boundary parallel annuli).  There are a number of possible
cases that arise in each step.  We will examine each in turn and
determine the effect on each of the constituents of $b$, denoting, for
example, that
$|S_+|$ goes up by at most two with the notation $|S_+| \up \leq 2$.

{\bf Step A: $\bdd$-reducing $W$ along $D$.}
\begin{enumerate}

\item  {\em $D$ is contained in a special component.}  Then $J \down 2$, $
|S_+| \up 1$, $ s \up 1$ so $b$ is unchanged.

\item  {\em $D$ is contained in a non-special component not in $W_-$.} 
Then $J \down 2$, $ |S_+| \up 1$, $ s \up 1$ so $b$ is unchanged.

\item  {\em $D$ lies in $W_-$.} Then $J \down 2$ and
$s \up 1$. The annulus abutting the parallelism created either abuts
$W_-$ on the other side (in which case $-m \up 1$) or does not (in
which case $-m \up 1$ and $-2n \down 2$).  Hence in any case $b$ does not
increase. 

\end{enumerate}

{\bf Step B: Inessential annuli $A$ removed.} There are various cases,
depending on what sort of component $C$ of $W - \Aaa$ lies on the other
side of the annulus $A$ (we will say that $A$ lies {\em on} the component
$C$) and whether or not the two boundary components of
$A$ lie on the same or different components of $C \cap \bdd_+ W$ (we will
say that $A$ is respectively non-separating or separating).  For example,
it may be that $C$ becomes special when $A$ is removed.  In this case $A$
is necessarily separating, since, by definition of ``special'', $\bdd C
\cap \bdd_+ W'$ is planar.  

\begin{enumerate}

\item  {\em $A$ is separating, doesn't lie on $W_-$, and isn't on a
component that becomes special.}  Then $|S_+| \down 2$, $ -2n \up 2$, $ s
\down 1$ (the last since the product component cut off by $A$ is special
by definition).  Hence $b$ drops by $1$.

\item  {\em $A$ lies on a component that becomes special.}  Then
$|S_+| \down 2$, $ -2n \up 2$ so $b$ is unchanged. 

\item  {\em $A$ is non-separating and doesn't lie on $W_-$.}  Then $|S_+|
\down 1$, $ -2n \up 2$, $ s \down 1$.  Hence $b$ is unchanged.

\item  {\em $A$ is separating and lies on $W_-$.}  Then $|S_+| \down 1$, $
-|S_-| \up 1$, $ -m \up 1$, $ s \down 1$.  Hence $b$ is unchanged.

\item  {\em $A$ is non-separating and lies on $W_-$.}  Then $|S_+| \down
1$, $ genus (S_-) \up 1$, $ -m \up 1$, $ s \down 1.$  Hence $b$ is
unchanged.

\end{enumerate}

Thus we conclude that $b' \leq b$ as required.  \qed

\begin{lem} \label{lem:count} With the hypotheses and notation of Lemma
\ref{lem:subcount}, let also $a$ denote the number of annuli
components of $S_+$, as in Lemma \ref{lem:generalcount}.  Then $\Aaa$ has
at most
$\frac{3}{4} J(W) +
\frac{1}{2} (a + s)$ components.
\end{lem}

\pf  By the special assumption that any orientable properly
embedded surface in $M^j$ has an even number of $\bdd$-components, each
component of $S_+$ is either an annulus or has Euler characteristic no
greater than $-2$.  Hence we have
\begin{equation}
 |S_+| \leq a - \frac{\chi(S_+)}{2} = a - \frac{\chi(\bdd_+ W)}{2} + 
\frac{\chi(S_-)}{2} 
\label{eq:S+}
\end{equation} On the other hand, $$- |S_-| + genus(S_-) + m + p = -
\frac{\chi(S_-)}{2}$$ so
\begin{equation}
 - |S_-| + genus(S_-) - m + p + \frac{\chi(S_-)}{2} = - 2m.
\label{eq:S-}
\end{equation} Adding inequality \ref{eq:S+}, equality \ref{eq:S-}, and
the equality
$$\frac{\chi(\bdd_- W)}{2} - \frac{\chi(S_-)}{2} = 
\frac{J}{2} + \frac{\chi(\bdd_+ W)}{2} - \frac{\chi(S_-)}{2}$$ we get
$$|S_+| - |S_-| + genus (S_-) - m + p + \frac{\chi(\bdd_- W)}{2} \leq a +
\frac{J}{2} - 2m$$ so from Lemma \ref{lem:subcount} we get
$$\frac{3J}{2} + a - 2m - 2n + s \geq 0$$ as required. 

\qed

\begin{lem} \label{lem:supercount} With the hypotheses and notation of
Lemma \ref{lem:count}, let also $I$ be the sum of the indices of the 
special components of $W -
\Aaa$ (See Definition \ref{def:special}). Then $\Aaa$ has at most
$\frac{3}{4} J(W) +
\frac{1}{2} (a + I)$ components.
\end{lem}

\pf Let $c = \frac{3}{2} J(W) + a + I - 2|\Aaa|$; it suffices to
show that $c \geq 0$. The proof is by induction on
$J$.  If there are no special components (e. g. when $J = 0$), the result
follows immediately from Lemma \ref{lem:count}.  Otherwise, let $D$ be a
$\bdd$-reducing disk for $W$ that lies in a special component $U$ of $W -
\Aaa$.  We will show that $\bdd$-reduction of $W$ along $D$ and the
removal of any resulting inessential annuli cannot increase 
$c$ in the compression body (or compression bodies) that arise.  

The first step, cutting $W$ along $D$, decreases $\frac{3}{2} J(W)$ by
$3$, raises $I$ by $3$ and raises $a$ by the number $\pi = 1$ or $2$ of
parallelisms in the two components $U - D$.  Thus $c - \pi$ is unchanged. 

Now examine the result of removing an inessential annulus $A \in \Aaa$
that cuts off a parallelism, an annulus such as is created when $W$ is
cut along $D$.  Let $U'$ be the component of $W - \Aaa$ that abuts $A$ on
the opposite side from $U$.  $U'$ may become special when $A$ is
removed.  There are several cases:
\begin{enumerate}

\item {\em $A$ was not adjacent to an annulus component of $S_+ \subset
\bdd W_+$}  Then even if $U'$ becomes special, $I(U') \leq 0$ so $-2|\Aaa|
\up 2$,   $a \down 1$, $I \down \geq 2$ (since the index of a parallelism
is $2$), $-\pi \up 1$.  Hence $c - \pi \down \geq 0$.

\item {\em $A$ is adjacent to an annulus component of $S_+ \subset
\bdd W_+$, $U'$ is not special or $I(U') \leq 0$.}  Then $-2|\Aaa|
\up 2$,   $a \down 2$, $I \down \geq 2$, $-\pi \up 1$.  Hence $c - \pi
\down \geq 1$.

\item {\em $A$ is adjacent to an annulus component of $S_+ \subset
\bdd W_+$, $U'$ is special and $I(U') = 1$.}   Then $-2|\Aaa|
\up 2$,   $a \down 2$, $I \down 1$, $-\pi \up 1$.  Hence $c - \pi$ is
unchanged.

\item {\em $A$ is adjacent to an annulus component of $S_+ \subset
\bdd W_+$, $U'$ is special and $I(U') = 2$.}  Then $U'$ becomes a
parallelism and $-2|\Aaa| \up 2$,   $a \down 2$, both $I$ and $-\pi$
unchanged.  Hence $c - \pi$ is unchanged.

\end{enumerate}

So $c - \pi$ never increases and may decrease.  Continue removing annuli
cutting off parallelisms until $\pi = 0$, as it was before the disk $D$
was removed.  In the end, $c$ also will not have increased.  \qed
 
\bigskip

Return now to the original context:  $M$ is a compact orientable
$3$-manifold and $\Aaa$ is a  complete properly imbedded
collection of essential annuli in $M$ satisfying the parity condition
(see Definition \ref{def:paritycond}).   Suppose
$$(V_1 \cup_{S_1} W_1)
\cup_{F_1} \dots \cup_{F_{m-1}} (V_m \cup_{S_m} W_m)$$ is a strongly
irreducible generalized Heegaard splitting of $M$ isotoped so that $(\Fff
\cup \Sss)
\cap A$ consists only of curves essential in both $\Fff
\cup \Sss$ and $\Aaa$, and such that this number is minimal.  Let
$M^1,...,M^n$ be the components into which $M$ is divided by the
family of annuli $\Aaa$. We will soon apply Lemma
\ref{lem:supercount} to the annuli $\Aaa \cap V_i$.  

Notice that the index of a special component of $V_i - \Aaa$
is non-positive unless the component $U$ is of index $1$, when $\bdd_S
U = \bdd U \cap S_i$ is a $4$-punctured sphere.  Since $\bdd_S U$ visibly
compresses in $U$, it follows from strong irreducibility that $\bdd_S U$
also compresses in $W_i$.  It then follows from a standard outermost arc
argument, that $\bdd_S U \cap S_i$ compresses in $W_i -
\Aaa$.  The result is two annuli.  If the two annuli are parallel, this
implies that $U$ is exceptional (see Definition \ref{defn:wide}).  If they
are not, then the union of the two annuli, together with the annuli $\bdd
U - \bdd_S U$, cuts off a
$2$-bridge knot or link complement from $M - \Aaa$, necessarily one of the
$M^j$ since $\Aaa$ is complete. 

\begin{defn} A component of $M^j$ which is a $2$-bridge knot or
link complement containing, as above, an index $1$ special component of
some
$V_i - \Aaa$ is called an \underline{interlaced} $2$-bridge complement.
\end{defn}

The following theorem now improves Theorem \ref{thm:generalcount}, when
$\Aaa \subset M$ satisfies the parity condition:

\begin{thm} \label{thm:count}
For $M$ and  $\Aaa$ as described above, let $k$ be the number of 
dipping annuli in $\Aaa - \Sss$, $w$ be the number of wide annuli in
$\Sss - \Aaa$, and
$e_V$ be the number of exceptional components cut off by $\Aaa$ from all
the $V_i$ (see Definitions
\ref{defn:dipping}, \ref{defn:wide}).  Let $c \leq m, n$ be the number of
interlaced $2$-bridge complements among the $M^j$.  Then
$$k \leq \frac{3}{4}(\Ss_{i=1}^m J(V_i)) + \frac{w + e_V + c}{2}.$$
\end{thm}

\pf  We will apply Lemma \ref{lem:supercount} to the annuli $\Aaa
\cap V_i$.  

Suppose, as a first simplification, that each
special component of each $V_i - \Aaa$ has non-positive index, so in
particular $e = c = 0$.  Then Lemma \ref{lem:supercount} says that the
number of dipping annuli in $\Aaa
\cap V_i$ is at most $\frac{3}{4}J(V_i) + \frac{a}{2}$, where $a$ counts
the number of annuli components of $S_i - \Aaa$ with a property (defined
in Lemma \ref{lem:count}) that assures that they are adjacent to dipping
annuli and therefore wide.  That is, when Lemma \ref{lem:supercount} is
applied to $V_i$ the inequality remains true if $a$ is replaced by the
number of wide annuli in $S_i$.  Summing over all $V_i$ we get
$k \leq \frac{3}{4}(\Ss_{i=1}^m J(V_i)) + \frac{w}{2}$ as required.  

Now examine what happens if $V_i - \Aaa$ has a special component $U$ of
positive index.  Then, as explained above, $U$ is of index $1$ and so is
either an exceptional component or part of an interlaced $2$-bridge
complement. In either case, the positive index contributed
by $U$ to $I$ in Lemma \ref{lem:supercount} is incorporated into $e_V +
c$. \qed

\begin{lem} \label{lem:partial}
For $M$ and  $\Aaa$ as described above, let $M^-$ be the manifold
obtained from $M$ by replacing each interlaced $2$-bridge 
complement $M^j$ with an unknot or unlink complement (matching
longitudes).  Then
$M^-$ has a (perhaps weakly reducible) generalized Heegaard splitting
$$(V^-_1
\cup_{S_1} W^-_1)
\cup_{F_1} \dots \cup_{F_{m-1}} (V^-_m \cup_{S_m} W^-_m)$$ so that each
$V^-_i \cong V_i$ and $W^-_i \cong W_i$.
\end{lem}

\pf  One way of constructing $M^-$ would be, on each index $1$
special component $U \subset (V_i - \Aaa)$ that comes from an interlaced
$2$-bridge complement, reglue $V_i$ to $W_i$ differently along
$\bdd_S U = \bdd U \cap S_i$.  Depending on how the $\bdd$-reducing
disk $D_W$ for $W_i$ separates pairs of boundary components of the
$4$-punctured sphere 
$\bdd_S U$, this regluing can be done so that the $\bdd$-reducing disk
$D_U$ of $U$ intersects $\bdd D_W$ in either two points or none. 
In the former case, the replacement is by an unknot complement and $U$
becomes exceptional (see Definition
\ref{defn:wide}).  In the latter, the replacement is by the complement
of an unlink of two components, and $M^-$ is reducible.
\qed

\begin{cor} \label{cor:special}
Let $M$ and  $\Aaa$ be as described above, with $\Aaa \subset M$
a complete collection of annuli satisfying the parity condition. Order the
summands so that the interlaced
$2$-bridge complements among the $M^j$, if any, are $M^{q+1}, \ldots
,M^n$. Then there is a Heegaard splitting
$V_-^j \cup_{S_-^j} W_-^j$ for each
$M^j$ so that
$$\frac{5}{2}(\Ss_{i=1}^m J(V_i)) + (n-q) \geq \sum_j J(V_-^j).$$
\end{cor}

\pf  We may assume that the surfaces $\Fff \cup \Sss$ from the
generalized Heegaard splitting and the annuli
$\Aaa$ have been isotoped so that $(\Fff \cup \Sss) \cap \Aaa$ consists
only of essential curves and a minimal number of them. Let $k$ denote the
number of dipping annuli in $\Aaa \cap (\cup_i V_i)$, $w$ denote the
number of wide annuli in $\Sss - \Aaa$, $e$ denote the total number of
exceptional components in the $V_i - \Aaa$ and the $W_i - \Aaa$, and
$e_V$ denote only the number that lie in the $V_i -
\Aaa$.  The construction in Theorem \ref{thm:bif} yields generalized
Heegaard splittings for $C(K^j)$ such that $2k + \sum J(V_i) \geq
\sum_{i,j} J(V_i^j)$.  

Now substitute for $k$ from Theorem \ref{thm:count} to get
$$\frac{5}{2}(\Ss_{i=1}^m J(V_i)) + w + e_V + (n-q) \geq \sum_{i,j}
J(V_i^j).$$  According to Lemma \ref{lem:dest} the induced 
Heegaard splittings for the $C(K^j)$ can be destabilized at
least $\frac{w}{2} + e$ times, yielding Heegaard splittings $V_-^j
\cup_{S_-^j} W_-^j$ for each $C(K^j)$ with $$\sum_j J(V_-^j) \leq
\sum_{i,j} J(V_i^j) - w - 2e \leq \sum_{i,j} J(V_i^j) - w - e_V.$$ It
follows that 
$$\frac{5}{2}(\Ss_{i=1}^m J(V_i)) + (n-q) \geq \sum_j J(V_-^j)$$ as
required. \qed

\section{Tunnel Numbers}
\label{sec:tun}

In this section we apply the results above to tunnel numbers of composite
knots.  The results can also be formulated for composite links, though
the statements are sometimes a bit more cumbersome.

%\begin{thm} \label{thm:sum}  Let $K^1, \ldots, K^n$ be prime knots, and
%suppose
%$K^{p+1}, \ldots, K^n$ are the $2$-bridge knots among the $K^j$.  Then
%there is some
%$q$ such that $p \leq q \leq n$ and
%\[ t(K^1 \# \ldots \# K^n) \left\{ 
%\begin{array}{ll}
%\geq \frac{2}{5}(t(K^1) + \ldots + t(K^q)) + \frac{n-q}{5} & and \\
%\geq (n-q) & and \\
%\geq t(K^1 \# \ldots \# K^q) &
%\end{array} \right. \]
%\end{thm}

\begin{thm} \label{thm:sum}  Let $K^1, \ldots, K^n$ be prime knots, and
suppose
$K^{p+1}, \ldots, K^n$ are the $2$-bridge knots among the $K^j$.  Then
there is some
$q$ such that $p \leq q \leq n$ and $t(K^1 \# \ldots \# K^n)$ is no
smaller than any of:

\begin{itemize}
\item  $\frac{2}{5}(t(K^1) + \ldots + t(K^q)) + \frac{n-q}{5}$
\item  $(n-q)$
\item  $t(K^1 \# \ldots \# K^q)$
\end{itemize}
\end{thm}

\pf By Lemma \ref{lem:ind} there is a generalized Heegaard
splitting $$C(K^1 \# \ldots \# K^n) \cong (V_1 \cup_{S_1} W_1) \cup_{F_1}
\dots \cup_{F_{m-1}}  (V_m \cup_{S_m} W_m)$$ for which $\sum_{i=1}^m
J(V_i) = 2t(K^1 \# \ldots \# K^n)$. Let $\Aaa$ be the $n-1$ annuli that
decompose
$C(K^1 \# \ldots \# K^n)$ into the complements of the constituent prime
knots. Clearly the non-trivial element of $H^1(C(K^1 \# \ldots \# K^n),
{\mathbb Z}_2)$ restricts to the parity class on meridians, so the parity
condition holds for $\Aaa \subset C(K^1 \# \ldots \# K^n)$.  

Let $q \geq p$ be the number of components of the
$C(K^j)$ that are not interlaced $2$-bridge complements.
The family $\Aaa$ is complete since each $K^j$ is prime. It
then follows from Corollary \ref{cor:special} that 
$$\frac{5}{2}(\Ss_{i=1}^m J(V_i)) + (n-q) \geq \sum_j J(V_-^j).$$  Since
$2$-bridge knots have tunnel number one, we further know that, for $j
\geq p$, hence $j \geq q$, that $J(V_-^j) \geq 2$.  Hence another way to
write the inequality is $$\frac{5}{2}(\Ss_{i=1}^m J(V_i)) \geq
\sum_{j=1}^q J(V_-^j) + (n-q)$$ or, substituting from Lemma \ref{lem:ind}
and dividing by two, $$\frac{5}{2}t(K^1 \# \ldots \# K^n) \geq
t(K^1) + \ldots + t(K^q) + \frac{n-q}{2}.$$  This gives the first
inequality.

The second inequality follows from the fact that each $V_i$ contributes
at least $2$ to $\Ss_{i=1}^m J(V_i)$ and the number $n-q$ of interlaced
$2$-bridge knot complements can be no bigger than $m$.  

The last inequality follows from Lemma \ref{lem:partial}. Here $M^-$ is
the manifold obtained by replacing each $C(K^j), q+1 \leq j \leq n$ with
the complement of the unknot.  Hence $M^- = C(K^1 \# K^2 \# \ldots
K^q)$.  It has a generalized Heegaard splitting whose constituent
compression bodies are the same as that of the untelescoped minimal genus
Heegaard splitting for $M$ itself, so when the splitting of $M^-$ is
amalgamated, the result is a Heegaard splitting of genus no higher than
that of $M$.  
\qed

\begin{cor} Let $K^1, \ldots, K^n$ be prime knots.

\begin{enumerate}  

\item $t(K^1 \# \ldots \# K^n)  \geq \frac{1}{3}(t(K^1) + \ldots +
t(K^n))$.

\item If none of the $K^j$ are $2$-bridge knots then $t(K^1 \# \ldots \#
K^n) \geq \frac{2}{5}(t(K^1) + \ldots + t(K^n))$. 

\item $t(K^1 \# K^2) \geq \frac{2}{5}(t(K^1) + t(K^2))$.

\end{enumerate}
\end{cor}

\pf   The first statement of course is just Corollary
\ref{cor:generaltunnel} but here it can be seen (somewhat intriguingly) to
follow also from the first and second inequalities of Theorem
\ref{thm:sum} and the fact that, for any
$a, b \geq 0$,  $max\{ (\frac{2}{5}a + \frac{1}{5}b), b \} \geq
\frac{1}{3}(a+b)$.

The second statement follows immediately from Theorem
\ref{thm:sum} by setting $p = n$.  

The last is immediate if both or neither of the $K^j$ are $2$-bridge.  If
exactly $K^2$ is $2$-bridge, the result follows from the last inequality
of Theorem \ref{thm:sum}.
\qed

\newpage

\appendix

\renewcommand{\epsilon}{\varepsilon}
\renewcommand{\d}{\partial}
\newcommand{\half}{\frac{1}{2}}
\newcommand{\ol}{\overline}
\newcommand{\Int}{\rm Int}
\newcommand{\ds}{\displaystyle}
\newcommand{\vs}[1]{\vspace{#1in}}
\newcommand{\RR}{\mathbb R}

\font\blah=line10
\def\setdiff{\hbox{\blah \hskip 4pt \char'121\hskip 4pt}}

\section{\scshape Examples with many dipping annuli}
\begin{center}
ANDREW CASSON
\end{center}

\begin{thm}	\label{example}
For any number $c < 2$,
there is a 3-manifold
$$M = X_1 \cup Y_1 \cup X_2 \cup Y_2 \cup \dots \cup X_k \cup Y_k$$
where $X_i$ and $Y_i$ are compression bodies,
and an essential annulus $A$ properly embedded in $M$,
such that the following conditions are satisfied.

\begin{enumerate}

\item
$\d M = \d_-X_1$ is an incompressible torus and $\d_-Y_k = \emptyset$.

\item
$X_i \cap Y_i = \d_+X_i = \d_+Y_i$,
and $X_i \cup Y_i$ is strongly irreducible.

\item
$Y_i \cap X_{i + 1} = \d_-Y_i = \d_-X_{i + 1}$ is incompressible in $M$.

\item
The components of $A \cap X_i$, $A \cap Y_i$
are essential sub-annuli of $A$.

\item
If $X_i$ has $n_i$ 2-handles,
and $A \cap X_i$ contains $d_i$ components
which are disjoint from $\d_-X_i$,
then ${\ds\sum_{i = 1}^k d_i \geq c\sum_{i = 1}^k n_i}$.

\end{enumerate}

\end{thm}

Most of the pairs $(X_i, A \cap X_i)$ and  $(Y_i, A \cap Y_i)$
are homeomorphic to a standard model $(V_n, F_n)$,
constructed in the following lemma.

\begin{lem}	\label{unit}
For each integer $n > 0$ there is a compression body $V_n$,
a homeomorphism $h : \d_+V_n \rightarrow \d_+V_n$,
and a surface $F_n$ properly embedded in $V_n$,
satisfying the following conditions.

\begin{enumerate}

\item
$V_n$ has $n$ 2-handles and $\d_-V_n$ is connected.

\item
$V_n \cup_h V_n$ is strongly irreducible.

\item
$F_n$ consists of $2n - 1$ essential,
non-parallel annuli disjoint from $\d_-V_n$ 
together with $4n - 2$ essential, non-parallel annuli
meeting both $\d_+V_n$ and $\d_-V_n$.

\item
$h(F_n \cap \d_+V_n) = F_n \cap \d_+V_n$.

\item
$F_n \cup_h F_n$ consists of $4n - 2$ annuli,
each meeting both boundary components of $V_n \cup_h V_n$.

\end{enumerate}

\end{lem}

\pf 
Let $S^2$ be the unit sphere $\{(x, y, z) \in \RR^3: x^2 + y^2 + z^2 = 1\}$.
Let $r_x$, $r_y$, $r_+$ and $r_-$ be the reflections of $S^2$
in the planes $y = 0$, $x = 0$, $x = y$ and $x = -y$ respectively,
generating a group $\Gamma$ of order 8.
The fixed point sets of $r_x$ and $r_y$ are the
circles $C_x$ and $C_y$ where $S^2$ intersects the
$x z$ and $y z$ planes respectively.

Let $Q^+ = \{(x, y, z) \in S^2: x \geq 0, y \geq 0\}$.
Let $D_1, D_2, \dots D_n$ be disjoint closed disks
in the interior of $Q^+$ such that $r_+(D_i) = D_i$,
for example, small round disks in $S^2$
with centers evenly spaced on the semicircle
where $Q^+$ meets the plane $x = y$.
Let $D^+ = D_1 \cup D_2 \cup \dots \cup D_n$
and $D = D^+ \cup r_x(D^+) \cup r_y(D^+) \cup r_-(D^+)$;
then $D$ is a $\Gamma$-invariant disjoint union of $4n$ disks in $S^2$.
(See Figure 10.) Set $P^+ = \ol{Q^+ \setdiff D^+}$ and $P = \ol{Q \setdiff
D}$, so $P$ is a $4n$-punctured sphere on which $\Gamma$ acts.
Let $S = P \cup_{\d} P'$,
where $P'$ is a connected surface with $\d P' = \d P$
such that the action of $\Gamma$ extends over $P'$.
For example, one could choose $P' = P$;
then $S$ is the double of $P$, a surface of genus $4n - 1$.

\begin{figure}
\centering
\includegraphics[width=.6\textwidth]{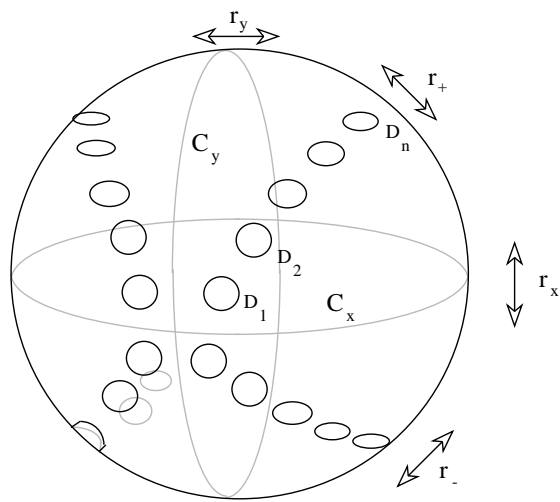}
\caption{The $4n$ disks $D$, with $n = 5$.}
\end{figure} 

Let $C^+_x = P^+ \cap C_x$, and let
$R = (P^+ \times 0) \cup (C^+_x \times I) \cup (P^+ \times 1)$
as a subset of of $\d(P^+ \times I)$.
There is a homeomorphism $g : R \rightarrow P^+ \cup r_x(P^+)$
such that $g(p, 0) = p$ and $g(p, 1) = r_x(p)$
for all $p \in P^+$ outside a small neighborhood $N$ of $C^+_x$,
and  $r_x g(p, t) = g(p, 1 - t)$
for every point $(p, t) \in R$.
Let $V_n$ be the 3-manifold obtained from the disjoint union
of $S \times I$ and $P^+ \times I$
by identifying each point $(p, t) \in R \subset P^+ \times I$
with $(g(p, t), 1) \in S \times 1 \subset S \times I$.
By using a collar neighborhood of $R$ in $P^+ \times I$,
construct an embedding $e : P^+ \times I \rightarrow V_n$
such that
$e(p, t) = (g(p, t), 0) \in S \times 0$ for all $(p, t) \in R$,
and $e(P^+ \times I)$ is disjoint from the subset
$(S \setdiff (P^+ \cup r_x(P^+)) \times I$ of $V_n$.

Let $I_1, I_2, \dots I_n$ be disjoint arcs properly embedded in $P^+$,
with end-points on $C^+_x$ and such that
$I_j$ separates $D_j$ from each $D_i$ with $i \neq j$.
Then $e(I_j \times I)$
is an embedded disk in $V_n$ with boundary in $S \times 0$.
This exhibits $V_n$ as a compression body
with $\d_+V_n = S \times 0$ and with
$e(I_1 \times I), e(I_2 \times I), \dots e(I_n \times I)$
as the cores of the 2-handles. (See Figure 11.)
Since $P'$ was chosen to be connected,
$\d_-V_n$ is connected.

\begin{figure}
\centering
\includegraphics[width=.6\textwidth]{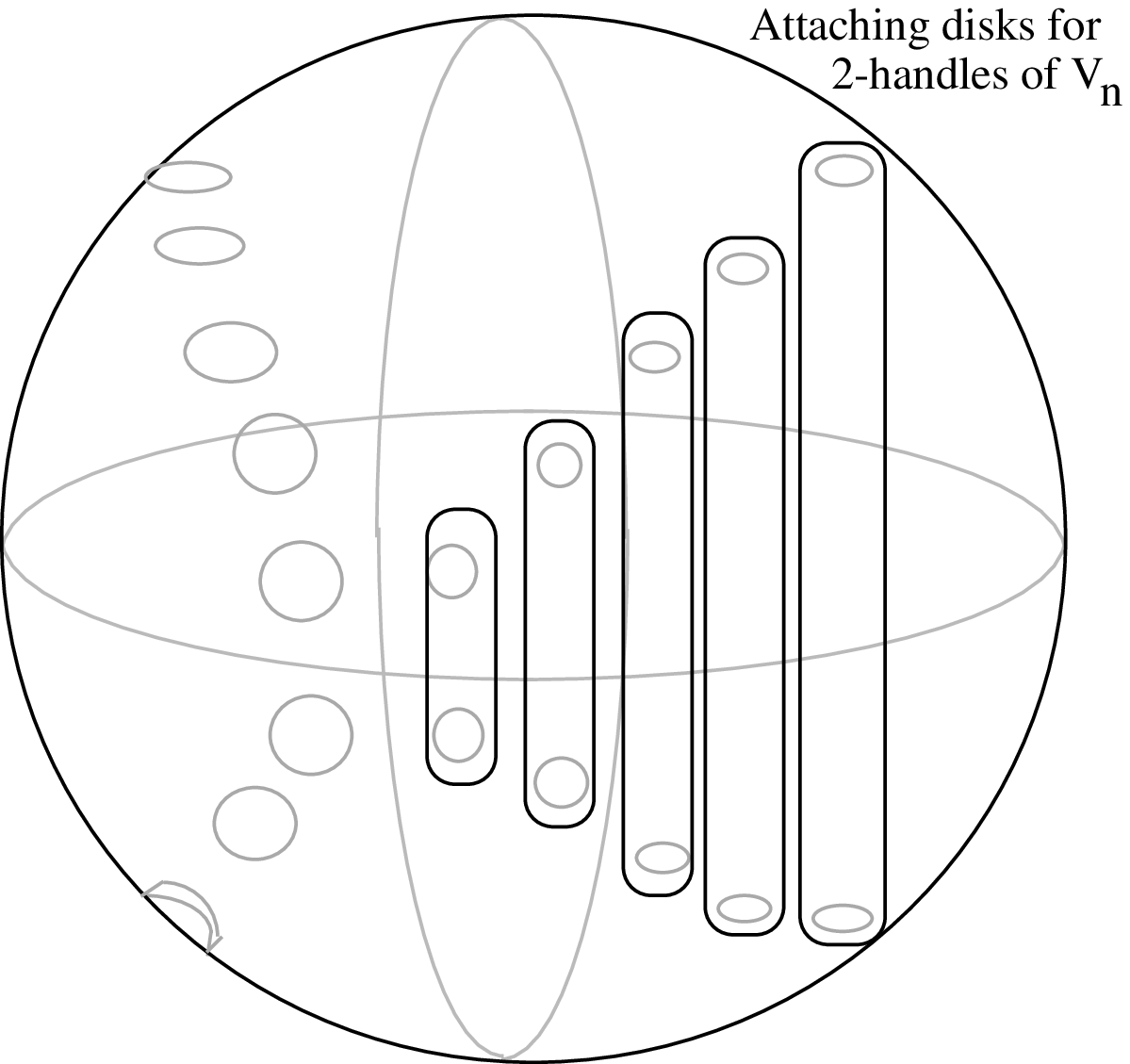}
\caption{The attaching disks $e(\bdd(I_j \times I))$}
\end{figure} 

For $1 \leq i \leq n$ let $C_i$ be a simple loop in $P^+$
that is parallel to the boundary component $\d D_i$
and invariant under $r_+$.
Enlarge $C_1, C_2, \dots C_n$ to
a maximal collection $C_1, C_2, \dots C_{2n - 1}$
of disjoint essential non-parallel simple loops in $P^+$
in such a way that $r_+(C_i) = C_i$ for all $i$. (See Figure 12.)
Set $A_i = e(C_i \times I)$;
then $A_1, A_2, \dots A_{2n - 1}$
are disjoint properly embedded annuli in $V_n$.
If $C_i$ is chosen disjoint from the neighborhood $N$ of $C^+_x$,
then $\d A_i = C_i \cup r_x(C_i)$.
$V_n$ also contains disjoint properly embedded annuli
$A'_1, A'_2, \dots A'_{2n - 1}$ and $A''_1, A''_2, \dots A''_{2n - 1}$,
where $A'_i = r_y(C_i) \times I$ and $A''_i = r_-(C_i) \times I$.
Let $F_n$ be the union of the annuli
$A_i$, $A'_i$ and $A''_i$ for $1 \leq i \leq 2n - 1$.

\begin{figure}
\centering
\includegraphics[width=.6\textwidth]{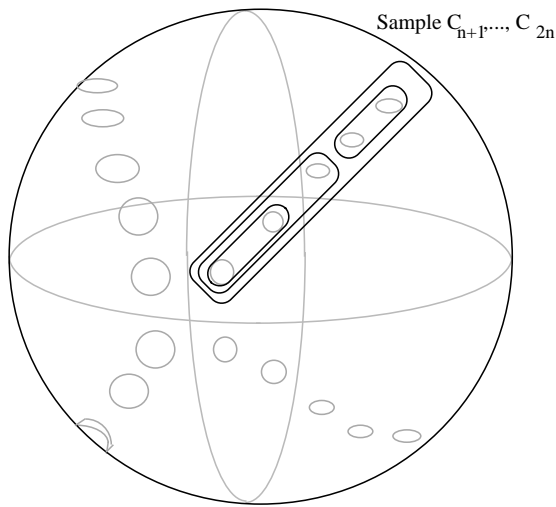}
\caption{A maximal collection of of disjoint essential non-parallel loops}
\end{figure} 

Let $h = t_x t_y r_+ : S \rightarrow S$,
where $t_x$ and $t_y$ are Dehn twists in $C_x$ and $C_y$ respectively.
By Corollary~\ref{strong} below, $V_n \cup_h V_n$ is strongly irreducible.
The intersection of $F_n$ with $\d_+V_n$ is the union of the curves
$C_i$, $r_x(C_i)$, $r_y(C_i)$ and $r_-(C_i)$ for $1 \leq i \leq 2n - 1$,
which is invariant under $h$.
In $V_n \cup_h V_n$,
each of $(A_i \cup A'_i) \cup_h (A_i \cup A'_i)$ and $A''_i \cup_h A''_i$
is a single properly embedded annulus.
\qed

%%%%%%%%%%%%%%%%%%%%%%%%%%%%%%%%%%%%%%%%%%%%%%%%%%%%%%%%%%%%%%%%%%%%%%%%%%%%%%%
%%%%%%%%%% ALSO NOTE THAT THE COMPONENTS OF F_n \cup_h F_n INDUCE THE %%%%%%%%%
%%%%%%%%%% IDENTITY PERMUTATION ON THE COMPONENTS OF F_n \cap \d_-V_n %%%%%%%%%
%%%%%%%%%%%%%%%%%%%%%%%%%%%%%%%%%%%%%%%%%%%%%%%%%%%%%%%%%%%%%%%%%%%%%%%%%%%%%%%

\begin{lem}	\label{wave}
Let $N$ be a regular neighborhood of $C_x \cup C_y$ in $S$.
Suppose $L$ is an essential simple loop in $S$ which
bounds a disk in $V_n$ and intersects $\d N \cup C_y$ minimally.
Then some component of $L \cap N$
intersects $C_x$ and is disjoint from $C_y$.
\end{lem}

\pf
A similar argument is given in \cite[Appendix]{MS}.
It follows from the construction of $V_n$ as
$S \times I \cup_{g} P^+ \times I$
that $\pi_1(V_n)$ is an HNN extension of $\pi_1(S \setdiff C_x)$,
so $\pi_1(S \setdiff C_x)$ injects into $\pi_1(V_n)$.
Moreover, every loop in $S$ which bounds a disk in $V_n$
has zero homological intersection number with $C_x$.
Therefore every essential loop in $S$ which bounds a disk in $V_n$
intersects $C_x$ at least twice.
Observe that $C_y$ bounds a disk in $V_n$
and intersects $C_x$ in exactly two points;
the regular neighborhood $N$ of $C_x \cup C_y$ is a 4-punctured sphere.

Suppose that $L$ is an essential simple loop in $S$ which
bounds a disk in $V_n$ and intersects $\d N$ minimally,
but every component of $L \cap N$
either intersects $C_y$ or is disjoint from $C_x$.
Since some component of $L \cap N$ intersects $C_x$,
$L$ must intersect $C_y$ also.

Therefore $L$ contains an arc $\alpha$ (a wave)
with $\alpha \cap C_y = \d\alpha$ and such that
if $\beta$ is either arc of $C_y \setdiff \d\alpha$
then $\alpha \cup \beta$ bounds a disk in $V_n$.
Each component of $\alpha \cap N$ not containing an end-point of $\alpha$
is a component of $L \cap N$ which is disjoint from $C_y$,
and is therefore disjoint from $C_x$ also.
Therefore every component of $\alpha \cap N$
contains an end-point of $\alpha$.

It follows that $(\alpha, \d\alpha)$ is isotopic in $(S, C_y)$
to an arc $\alpha'$ not meeting $C_x$;
let $\beta'$, $\beta''$ be the arcs of $C_y \setdiff \d\alpha'$.
Then $\alpha' \cup \beta'$, $\alpha' \cup \beta''$ are
essential simple closed curves in $S$ bounding disks in $V_n$,
and at least one intersects $C_x$ in less than two points,
a contradiction.
\qed

\begin{cor}	\label{strong}
With $V_n$ and $h$ as in Lemma~\ref{unit},
$V_n \cup_h V_n$ is strongly irreducible.
\end{cor}

\begin{lem}	\label{end}
Let $S$ be a closed orientable surface,
possibly disconnected,
and let $C_1, C'_1, \dots C_k, C'_k$ be
disjoint essential simple closed curves on $S$.
Then there is a compact 3-manifold $N$
with incompressible boundary $\d N = S$,
containing disjoint annuli
$A_1, A_2, \dots A_k$
such that $\d A_i = C_i \cup C'_i$.
\end{lem}

\pf
Let $A_1, A_2, \dots A_k$ be annuli, and construct
$$N' = (S \times I) \cup (A_1 \times I) \cup (A_2 \times I) \cup
\dots \cup (A_k \times I)$$
where $(A_i \times I) \cap (S \times I) = (\d A_i) \times I$
is a regular neighborhood of $C_i \cup C'_i$ in $S \times 1$.
The boundary of $N'$ consists of $S$ together with a closed surface $S'$;
since $C_i$ and $C'_i$ are essential, $S$ and $S'$ are incompressible in $N'$.
Choose a a compact 3-manifold $N''$
with incompressible boundary $\d N'' = S'$,
and set $N = N' \cup_{S'} N''$.
\qed

\bigskip

\noindent{\bf Proof of Theorem \ref{example}}
Choose $n > (1 - c / 2)^{-1}$
and let $V_n$, $h$ and $F_n$ be as in Lemma~\ref{unit}.
Let $C$, $C'$ be parallel essential simple closed curves on a torus $T$.
Let $S$ be the disjoint union of $T$ and $\d_-V_n$,
and let $C_1, C'_1, \dots C_{2n}, C'_{2n}$
be the components of $(F_n \cap \d_-V_n) \cup C \cup C'$,
in any order.
By Lemma~\ref{end},
there is a 3-manifold $N$ with incompressible boundary $T \cup \d_-V_n$,
containing disjoint annuli $A_1, A_2, \dots A_{2n}$ such that
$$(F_n \cap \d_-V_n) \cup C \cup C'
= \d A_1 \cup \d A_2 \cup \dots \cup \d A_{2n}.$$

If instead $S = \d_-V_n$ and
$C_1, C'_1, \dots C_{2n - 1}, C'_{2n - 1}$
are the components of $F_n \cap \d_-V_n$ in any order
then Lemma~\ref{end}
gives a 3-manifold $N'$ with incompressible boundary $\d_-V_n$
and containing disjoint annuli $A'_1, A'_2, \dots A'_{2n - 1}$ such that
$$F_n \cap \d_-V_n = \d A'_1 \cup \d A'_2 \cup \dots \cup \d A'_{2n - 1}.$$

Choose the orderings in such a way that
if $F = A_1 \cup A_2 \dots \cup A_{2n}$
and $F' = A'_1 \cup A'_2 \cup \dots \cup A'_{2n - 1}$,
then $F \cup F' \subset N \cup_S N'$ is a single annulus.

By \cite[Lemma 3]{SS1},
$(N; T, \d_-V_n)$ has a strongly irreducible untelescoped Heegaard splitting
$N = X_1 \cup Y_1 \cup \dots \cup X_k \cup Y_k$
such that the components of $A_i \cap X_j$ and $A_i \cap Y_j$
are essential sub-annuli of $A_i$.
Similarly,
$(N'; \d_-V_n, \emptyset)$ has an untelescoped Heegaard splitting
$N' = X'_1 \cup Y'_1 \cup \dots \cup X'_{k'} \cup Y'_{k'}$
with similar properties.

Let $q$ be the total number of 1-handles in
these generalized Heegaard splittings for $N$ and $N'$.
Choose $m > 2q$ and let $m(V_n \cup_h V_n)$ denote the union
$(V_n \cup_h V_n) \cup (V_n \cup_h V_n) \cup \dots \cup (V_n \cup_h V_n)$
of $m$ copies of $V_n \cup_h V_n$.
Set $M = N \cup m(V_n \cup V_n) \cup N'$
with the generalized Heegaard splitting
$$X_1 \cup Y_1 \cup \dots \cup X_k \cup Y_k \cup m(V_n \cup_h V_n)
\cup X'_1 \cup Y'_1 \cup \dots \cup X'_{k'} \cup Y'_{k'}$$
and set $A = F \cup m(F_n \cup F_n) \cup F'$,
a single annulus properly embedded in $M$.

The total number of 1-handles is $mn + q$,
and there are at least $(2n - 1)m$ ``dipping'' annuli.
Since $n \geq (1 - c/2)^{-1}$ and $m \geq 2q$,
$(2 - c)mn \geq 2m \geq m + 2q \geq m + cq$,
so $(2n - 1)m \geq c(mn + q)$, as required.
\qed


\begin{thebibliography}{99}

\bibitem[BnO]{BnO} F.~Bonahon, J.P.~Otal, {\em Scindements de Heegaard des
espaces lenticulaires}, Ann. scient. \'{E}c. Norm. Sup. $4^e$
s\'{e}rie {\bf 16} (1983) 451--466.

\bibitem[BZ]{BZ} G.~Burde, H.~Zieschang, {\em Knots}
de Gruyter, Studies in Mathematics {\bf 5}, Berlin, New York.

\bibitem[CG]{CG} A.~Casson and C.~McA.~Gordon, {\em Reducing Heegaard
splittings}, Topology  and its applications {\bf} 27 (1987), 
275-283.

\bibitem[Ga]{Ga} D. Gabai, {\em Foliations and the topology of
$3$-manifolds}, J. Diff. Geom. {\bf 18} (1983) 445-503.

\bibitem[H]{H} J.~Hempel, {\em $3$-manifolds}, Annals of Math.
Studies {\bf 86} (1976), Princeton University Press.

\bibitem[J]{J} W.~Jaco, {\em Lectures on three-manifold Topology}.
Regional Conference Series in Mathematics {\bf 43} (1981), Amer. Math.
Soc.

\bibitem[Ko]{Ko} T.~Kobayashi, {\em A construction of arbitrarily high
degeneration of tunnel numbers of knots under connected sum}, J. Knot
Theory Ramifications {\bf 3} (1994), no. 2, 179-186.

\bibitem[Kw]{Kw} H.-Z.~Kowng, {\em Straightening tori in Heegaard
splittings}, Ph. D. thesis, U. C. Santa Barbara, 1994.

\bibitem[M1]{M1} K. Morimoto, {\em On composite tunnel number one links},
Topology Appl., {\bf 59}(1994), no. 1, 59--71.

\bibitem[M2]{M2} K. Morimoto, {\em There are knots whose tunnel numbers go
down under connected sum}, Proc. Am. Math. Soc, {\bf 123} (1995), no. 11,
3527--3532

\bibitem[MS]{MS} Y. Moriah and J. Schultens, {\em Irreducible Heegaard
splittings of Seifert fibered spaces are either vertical or horizontal},
Topology {\bf 37} (1998), no 5, 1089-1112.

\bibitem[R]{R} D. Rolfsen, {\em Knots and links}, Publish or Perish, Inc.
Houston, Texas

\bibitem[Sh]{Sh} M.~Scharlemann, {\em Local detection of strongly
irreducible Heegaard splittings}, Topology and its Applications {\bf 90}
(1998), 135-147

\bibitem[ST1]{ST1} M.~Scharlemann, A.~Thompson,{\em Thin position for
$3$-manifolds}, AMS Contemporary Math. {\bf 164} (1994), 231--238

\bibitem[SS1]{SS1} M.~Scharlemann, J.~Schultens, {\em The tunnel number
of the sum of $n$ knots is at least $n$}, to appear in Topology.

\bibitem[SS2]{SS2} M.~Scharlemann, J.~Schultens, {\em Comparing JSJ and
Heegaard structures of orientable $3$-manifolds}, preprint.

\bibitem[Sc]{Sc} J. Schultens, {\em Additivity of tunnel number for small
knots}, preprint.

\end{thebibliography}
\end{document}